\documentclass[11pt,a4paper]{article}
\usepackage[margin=1in]{geometry}
\usepackage{authblk}

\usepackage{microtype}
\usepackage{graphicx}
\usepackage{subfigure}
\usepackage{booktabs}
\usepackage{enumitem}
\usepackage{xcolor}
\usepackage{tikz}
\usepackage[ruled,vlined,linesnumbered]{algorithm2e}

\usepackage{hyperref}
\usepackage{natbib}

\usepackage{amsmath, amssymb, mathtools, amsthm, bm}

\theoremstyle{plain}
\newtheorem{theorem}{Theorem}[section]
\newtheorem{proposition}[theorem]{Proposition}
\newtheorem{lemma}[theorem]{Lemma}

\theoremstyle{definition}

\newtheorem{assumption}[theorem]{Assumption}
\theoremstyle{remark}

\newcommand{\ABSTRACT}[1]{\begin{abstract}#1\end{abstract}}

\usepackage{amsmath}
\allowdisplaybreaks

\usepackage{amssymb}

\usepackage{bm}
\usepackage{hyperref}
\usepackage{url}
\usepackage{booktabs} 

\usepackage{multirow}

\usepackage{rotating}

\usepackage{tcolorbox}

\newcommand{\bb}{\bm{b}}

\newcommand{\bI}{\bm{I}}

\newcommand{\bX}{\bm{X}}

\newcommand{\bz}{\bm{z}}

\newcommand{\bs}{\bm{s}}

\newcommand{\by}{\bm{y}}

\newcommand{\balpha}{\bm{\alpha}}

\newcommand{\bbeta}{\bm{\beta}}
\newcommand{\bomega}{\bm{\omega}}

\newcommand{\brho}{\bm{\rho}}

\newcommand{\bgamma}{\bm{\gamma}}

\newcommand{\bmu}{\bm{\mu}}
\newcommand{\bnu}{\bm{\nu}}
\newcommand{\bpi}{\bm{\pi}}

\newcommand{\bzeta}{\bm{\zeta}}

\newcommand{\bbR}{\mathbb{R}}

\newcommand{\calB}{\mathcal{B}}

\newcommand{\calJ}{\mathcal{J}}

\newcommand{\calN}{\mathcal{N}}

\def\1{\bm{1}}

\DeclareMathAlphabet{\mathsfit}{\encodingdefault}{\sfdefault}{m}{sl}
\SetMathAlphabet{\mathsfit}{bold}{\encodingdefault}{\sfdefault}{bx}{n}

\newcommand{\R}{\mathbb{R}}

\def\vert#1{\lvert #1 \rvert}

\definecolor{predcolor}{gray}{0.95}
\definecolor{scorecolor}{gray}{0.95}
\definecolor{riskcolor}{gray}{0.95}
\definecolor{transparentcolor}{gray}{0.95}

\definecolor{ForestGreen}{rgb}{0.13, 0.55, 0.13} 
\usepackage{color} 

\DeclareMathOperator{\conv}{conv}
\DeclareMathOperator*{\TopSum}{TopSum}
\DeclareMathOperator{\prox}{prox}
\undef\st
\DeclareMathOperator{\st}{s.\!t.\!}

\DeclareMathOperator{\epi}{epi}
\DeclareMathOperator{\ri}{rint}
\DeclareMathOperator{\rint}{rint}
\DeclareMathOperator{\inter}{int}

\DeclareMathOperator{\dom}{dom}
\DeclareMathOperator{\dist}{dist}
\DeclareMathOperator{\init}{init}

\usetikzlibrary{arrows.meta}

\providecommand{\argmin}{\operatorname*{arg\,min}}
\providecommand{\prox}{\operatorname{prox}}

\newtheorem*{rep@theorem}{\rep@title}

\newtheorem*{rep@lemma}{\rep@title}

\newtheorem*{rep@proposition}{\rep@title}

\title{GPU-friendly and Linearly Convergent First-order Methods for Certifying Optimal $k$-sparse GLMs}

\author[$\ast,\dagger$]{Jiachang Liu}
\author[$\ddagger$]{Andrea Lodi}
\author[$\dagger$]{Soroosh Shafiee}
\affil[$\ast$]{Center for Data Science for Enterprise and Society, Cornell University, Ithaca, USA}
\affil[$\dagger$]{\mbox{School of Operations Research and Information Engineering, Cornell University, Ithaca, USA}}
\affil[$\ddagger$]{Jacobs Technion-Cornell Institute, Cornell Tech and Technion--IIT, New York, USA\\
\texttt{\{jiachang.liu, andrea.lodi, shafiee\}@cornell.edu}}

\date{}

\begin{document}
\maketitle
\addtocontents{toc}{\protect\setcounter{tocdepth}{-10}}

\ABSTRACT{
We investigate the problem of certifying optimality for sparse generalized linear models (GLMs), where sparsity is enforced through a cardinality constraint.
While Branch-and-Bound (BnB) frameworks can certify optimality using perspective relaxations, existing methods for solving these relaxations are computationally intensive, limiting their scalability.
To address this challenge, we reformulate the relaxations as composite optimization problems and develop a unified proximal framework that is both linearly convergent and computationally efficient.
Under specific geometric regularity conditions, our analysis links primal quadratic growth to dual quadratic decay, yielding error bounds that make the Fenchel duality gap a sharp proxy for progress towards the solution set.
This leads to a duality gap-based restart scheme that upgrades a broad class of sublinear proximal methods to provably linearly convergent methods, and applies beyond the sparse GLM setting.
For the implicit perspective regularizer, we further derive specialized routines to evaluate the regularizer and its proximal operator exactly in log-linear time, avoiding costly generic conic solvers.
The resulting iterations are dominated by matrix--vector multiplications, which enables GPU acceleration.
Experiments on synthetic and real-world datasets show orders-of-magnitude faster dual-bound computations and substantially improved BnB scalability on large instances.
}

\section{Introduction}
\label{sec:introduction}
Sparse generalized linear models (GLMs) are essential tools in machine learning (ML) and statistics, with broad applications in healthcare, finance, engineering, and science. 
In this paper, we study the global optimization of cardinality-constrained GLMs; a formal formulation is given in Section~\ref{sec:sparse_glm_problem_setup}. 
While many scalable approaches rely on convex surrogates (e.g., lasso) or other heuristics, such approximations can yield poor solutions in the presence of high-dimensional or highly correlated features~\citep{liu2022fast, liu2024okridge}.
This limitation is particularly problematic in high-stakes applications like healthcare, where accuracy, reliability, and interpretability are essential.
Therefore, we emphasize the pursuit of certifiably optimal solutions to the exact $\ell_0$-constrained formulation.
Such problems are NP-hard, and a standard approach is to solve them with mixed-integer programming using branch-and-bound (BnB), which requires computing valid lower bounds at every node of the search tree.
However, the standard big--M relaxation is often weak, leading to loose lower bounds and slow pruning.
To address this issue, the perspective relaxation~\citep{ceria1999convex, gunluk2010perspective, atamturk2020supermodularity} has become a popular mechanism for strengthening node relaxations and improving BnB performance.
The perspective relaxation yields a tighter lower bound than the standard big--M relaxation~\citep[Lemma~2.1]{liu2025scalable} and can often be expressed as a conic convex program.

However, solving the perspective relaxations at scale is challenging.
One standard approach is to use the interior-point method (IPM), but it does not scale well due to its reliance on solving linear systems, which has a cubic complexity in the number of variables and cannot be easily parallelized on modern hardware like GPUs.
Additionally, IPM cannot be effectively warm-started.
Due to these limitations, there has been growing interest in solving such relaxations using first-order methods~\citep{hazimeh2022sparse, liu2024okridge}, which are more scalable and can be easily warm-started.
However, in order to truly benefit from GPUs, these first-order methods must avoid expensive subroutines that cannot be parallelized.
It is unclear how to design such first-order methods that can fully leverage GPU parallelism for perspective relaxations.

Another drawback of first-order methods is their slow convergence.
From a theoretical perspective, there is a lack of understanding regarding the convergence rates of these methods, particularly for obtaining the \emph{safe lower bounds} required for pruning in BnB.
Since iterative methods produce only approximate solutions, the objective value at any iterate is an upper bound~on the minimum of the relaxation, which is not a valid lower bound for the original problem.
While techniques exist to derive safe lower bounds from inexact solutions, the convergence rate of these bounds has not been addressed.
This is critical as the speed at which accurate lower bounds are obtained directly dictates the efficiency of the BnB algorithm.
Ideally, one would like to achieve a linear rate, but it is unclear if such rates are attainable for these perspective relaxations using first-order methods.

Both the aforementioned computational challenges and the theoretical gap on convergence rates give rise to the central question of our paper, namely

\begin{tcolorbox}[colback=blue!5!white,colframe=gray!75!black,title=Central Research Question]
    \label{tcolorbox:central_research_question}

    \emph{Can we solve the perspective relaxations within branch-and-bound and obtain valid lower bounds by designing a GPU-friendly first-order method with provable linear convergence?}

\end{tcolorbox}

Our answer to this question is positive.
Below, we list the main technical contributions we have developed in order to answer this question affirmatively.
\begin{enumerate}[label=$\diamond$]
    \item \emph{Composite Reformulation:} We reformulate the perspective relaxation as an unconstrained convex composite problem, introducing a novel non-smooth regularizer implicitly defined by the perspective function.
    
    \item \emph{Geometric Analysis:} We analyze a general composite structure and uncover two favorable geometric properties.
    Namely, under certain regularity conditions, the primal problem satisfies the \emph{quadratic growth condition}, while the dual problem satisfies the \emph{quadratic decay condition}, a term we introduce as the dual counterpart to the primal quadratic growth.
    
    \item \emph{Linearly Convergent Restart Scheme:} Exploiting these geometric properties, we propose a novel \emph{gap-based restart scheme}.
    This generic framework, which applies broadly to convex composite problems beyond perspective relaxations, accelerates a broad class of proximal methods, including fixed-stepsize, linesearch-based, and linesearch-free variants, to achieve \emph{provable linear convergence rates} for both the primal and dual objectives and sequences.
    
    \item \emph{Efficient Implementation for Sparse GLMs:} We prove that the implicit perspective regularizer satisfies the required geometric conditions and show that both its function value and proximal operator can be evaluated \emph{exactly and efficiently} in log-linear time. This enables the use of efficient, GPU-friendly first-order methods without solving expensive conic subproblems.
    
    \item \emph{Empirical Performance:} We validate the practical efficiency of our approach on both synthetic and real-world datasets, showcasing substantial speedups of one to two orders of magnitude on CPUs and an additional order of magnitude on GPUs in computing dual bounds and certifying optimal solutions for large-scale sparse GLMs.
\end{enumerate}


\noindent\textbf{\emph{Relation to \citet{liu2025scalable}.}}
\citet{liu2025scalable} study proximal methods for the perspective relaxation at the root node.
This paper extends that conference version by studying perspective relaxations at \emph{every node}, adding the \emph{geometric analysis} and \emph{linearly convergent restart} scheme highlighted above, and reporting \emph{GPU-accelerated BnB} experiments for certifying optimality.


\subsection{Related Works}
\label{sec:related_work}

\paragraph{MIP for ML.}
Mixed integer programming (MIP) has been applied to
medical scoring systems~\citep{ustun2016supersparse, ustun2019learning, liu2022fasterrisk}, 
portfolio optimization \citep{bienstock1996computational,wei2022convex},
nonlinear identification systems~\citep{bertsimas2023learning, liu2024okridge},
decision trees~\citep{bertsimas2017optimal, hu2019optimal},
survival analysis~\citep{zhang2023optimal, liu2024fastsurvival},
hierarchical models~\citep{bertsimas2020sparse},
regression and classification models~\citep{atamturk2020safe, bertsimas2020sparse, bertsimas2020sparse1, bertsimas2020sparse2, hazimeh2020fast, xie2020scalable, atamturk2021sparse, dedieu2021learning, hazimeh2022sparse, liu2024okridge, guyard2024el0ps},
graphical models \citep{manzour2021integer, kucukyavuz2023consistent},
explainability~\citep{dash2018boolean, lodi2024one},
and outlier detection \citep{gomez2021outlier,gomez2023outlier}.
Most work targets high-quality feasible solutions, with few addressing optimality certification~\citep{tillmann2024cardinality}; we focus on computationally scalable optimality certification for sparse GLMs.

\paragraph{Perspective Formulations.} 
Applications of perspective relaxations like ours (Section~\ref{sec:sparse_glm_problem_setup}) date back to \citet{ceria1999convex}.
Perspective formulations for separable functions have been developed in \citep{gunluk2010perspective, xie2020scalable, wei2022ideal, bacci2019new, shafiee2024constrained}, and for rank-one functions in \citep{atamturk2020supermodularity, wei2020convexification, wei2022ideal, han2021compact, shafiee2024constrained}.
Our work uses separable perspective formulations arising from the Tikhonov regularization term in sparse GLMs.
We introduce a new implicit function based on the perspective function and derive its Fenchel conjugate, proximal operator, function evaluation, and geometric properties.

\paragraph{Lower Bound Calculation.}
Commercial MIP solvers typically iteratively linearize the objective function using the outer approximation method~\citep{kelley1960cutting} (via cutting planes) and solve the resulting linear programs~\citep{schrijver1998theory, wolsey2020integer}.
However, this approach often produces loose lower bounds, especially when high-quality cuts are not generated.
Alternatively, solvers may solve conic convex relaxations with IPMs~\citep{dikin1967iterative, renegar2001mathematical, nesterov1994interior}, but IPMs are not scalable.
Recent approaches use first-order methods, including subgradient descent~\citep{bertsimas2020sparse1}, alternating direction method of multipliers (ADMM)~\citep{liu2024okridge}, and coordinate descent~\citep{hazimeh2022sparse}.
Our work builds on this direction.
Our implicit formulation enables proximal gradient methods with low-cost subroutines that parallelize well on GPUs.
Crucially, we establish provable global linear convergence rates for derived safe lower bounds at every node in BnB, which previous first-order methods for this problem did not provide.



\paragraph{Proximal Gradient Methods}
For convex composite problems with a smooth term, proximal gradient descent (PGD) achieves a sublinear rate of $O(1/k)$~\citep{chen1997convergence,combettes2005signal}.
Applying Nesterov acceleration leads to the fast iterative shrinkage-thresholding algorithm (FISTA) with an improved sublinear rate of $O(1/k^2)$~\citep{nesterov1983method, beck2009fast, tseng2008accelerated}.
Recent parameter-free variants, notably the Adaptive Gradient Descent (AdProxGD)~\citep{malitsky2024adaptive} and the Adaptive Composite Fast Gradient Method (AC-FGM)~\citep{li2023simple}, remove manual step-size tuning.
Although these rates are optimal in general, first-order methods can converge \emph{linearly} under strict regularity conditions, such as strong convexity, the Polyak-{\L}ojasiewicz condition~\citep{karimi2016linear}, or the quadratic growth condition~\citep{luo1993error, drusvyatskiy2018error}.
Under the quadratic growth condition, restart schemes can accelerate momentum methods, ranging from heuristic schemes~\citep{o2015adaptive} to more systematic schemes~\citep{fercoq2019adaptive} based on proximal gradient mappings or function values~\citep{roulet2017sharpness}. 
We propose a restart scheme \emph{based on the duality gap}; unlike existing schemes that rely solely on primal information, it applies broadly and yields linear convergence.

\paragraph{Relationship to Dual-Based Strategies.}
Our use of the \emph{duality gap} for restarts shares conceptual roots with FDPG~\citep{beck2014fast}, but offers a new perspective and improves the convergence rate.
FDPG applies FISTA to the dual, yielding $\mathcal{O}(1/k^2)$ for the dual objective, but only $\mathcal{O}(1/k)$ for the recovered primal.
We proceed in the reverse direction by optimizing the primal and construct a dual sequence to monitor the duality gap.
With a gap-based restart strategy, we ensure \emph{both} the primal and dual sequences converge \emph{linearly}, outperforming the FDPG's sublinear rate.
Moreover, we do not assume the strong convexity, which makes the dual objective globally smooth.

\paragraph{GPU Acceleration.}
Recent work uses GPUs to accelerate continuous optimization problems, including linear programming~\citep{applegate2021practical, lu2023cupdlp}, quadratic programming~\citep{lu2023practical}, and semidefinite programming~\citep{han2024accelerating}.
A natural way to leverage GPUs for discrete problems is through GPU-based LPs within MIP solvers, as demonstrated by~\citet{de2024power} for clustering.
However, their approach may require iteratively generating many cutting planes to approximate the original objective.
In contrast, we develop a customized proximal gradient method that directly handles the nonlinear objective, while the computation is GPU-friendly since it mainly involves matrix--vector multiplication. 
Other first-order methods, such as ADMM~\citep{liu2024okridge} and coordinate descent~\citep{hazimeh2022sparse}, are less suitable for GPUs: ADMM requires solving linear systems, while coordinate descent is inherently sequential.

\section{Preliminaries and Problem Formulation}
In this section, we review the mathematical background needed to analyze the perspective relaxation of sparse GLMs within a BnB framework. We first introduce standard notation and convex analysis concepts used throughout the paper. We then state the sparse GLM problem and the perspective relaxation solved at each node of the BnB tree. We next recast the perspective relaxation as an \emph{unconstrained composite optimization problem}. Finally, we present its Fenchel dual, which forms the basis for our duality-gap-based restart scheme.

\subsection{Notation}
\label{sec:notation}
Throughout the paper, we work in the standard Euclidean space $\mathbb{R}^p$, where the inner product of two vectors $\bbeta, \bgamma \in \mathbb{R}^p$ is denoted by $\bbeta^\top \bgamma$, and the (Euclidean) norm of $\bbeta$ is denoted by $\|\bbeta\|_2$.
The extended real line is defined as $\overline{\mathbb{R}} := \mathbb{R} \cup \{ -\infty, +\infty \}$.
For any $p \in \mathbb N$, we define $[p] := \{ 1, \dots, p \}$. We also define $[p_1:p_2] := \{p_1, \dots, p_2\}$ for any $p_1, p_2 \in \mathbb N$ with $p_2 \geq p_1$.
The number of elements in a finite set $\mathcal{J}$ is denoted by $\left| \mathcal{J} \right|$.
For a set $\mathcal{J} \subseteq [p]$, the subvector $\bbeta_{\mathcal{J}}$ consists of the components of $\bbeta \in \mathbb{R}^p$ indexed by $\mathcal{J}$. For example, the subvector $\bbeta_{[p_1 : p_2]}$ consists of the components of $\bbeta$ from index $p_1$ to $p_2$.
We denote by $\bm{1}_p$ and $\bm{0}_p$ the vectors of all ones and all zeros in $\mathbb R^p$, respectively, and drop the subscript $p$ when the dimension is clear from context.
For a vector $\bbeta \in \mathbb R^p$, we use the notation $|\beta_{(1)}| \ge \cdots \ge |\beta_{(p)}|$ to denote the magnitudes of its entries sorted in nonincreasing order (ties broken arbitrarily).
Equivalently, $\beta_{(j)}$ denotes an entry of $\bbeta$ whose magnitude is the $j$-th largest among $\{|\beta_1|,\ldots,|\beta_p|\}$.
Given a set $\mathcal{C} \subseteq \mathbb{R}^p$, we denote by $\inter(\mathcal C)$ and $\rint(\mathcal C)$, the interior and the relative interior of $\mathcal C$, respectively. The distance function between $\bbeta \in \mathbb{R}^p$ and the set $\mathcal{C}$ is defined as $\dist(\bbeta, \mathcal{C}) := \inf_{\bgamma \in \mathcal{C}} \|\bbeta - \bgamma\|_2$.
The distance between two sets $\mathcal{B}, \mathcal{C} \subseteq \mathbb{R}^p$ is defined as $\dist(\mathcal B, \mathcal{C}) = \inf_{\bbeta \in \mathcal{B}} \dist(\bbeta, \mathcal{C})$.
We use $\delta_{\mathcal{C}}(\bbeta)$ to denote the indicator function of the set $\mathcal{C}$, defined as $\delta_{\mathcal{C}}(\bbeta) = 0$ if $\bbeta \in \mathcal{C}$; $= +\infty$ otherwise.

\subsection{Convex Analysis Background}
For an extended-valued function $h : \mathbb{R}^p \to \overline{\mathbb{R}}$, its (effective) domain is defined as $ \dom( h ) := \{ \bbeta \in \mathbb{R}^p \mid h ( \bbeta ) < +\infty \} $. The \emph{epigraph} of $h$ is defined as the set $\epi(h) = \{(\bbeta, t) \mid h(\bbeta) \leq t \}$.
The function $h$ is \emph{proper} if it never takes $-\infty$ and has nonempty domain. 
It is called \emph{closed} if its epigraph is closed. 
The \textit{Fenchel conjugate} of $h$ is defined as 
\begin{align*}
    h^* ( \balpha ) := \sup_{\bbeta \in \mathbb{R}^p} \ \bbeta^\top \balpha - h ( \bbeta ),
\end{align*}
and the \textit{proximal operator} of $h$ is defined as 
\[ \prox_h ( \bbeta ) := \argmin_{\bgamma \in \mathbb{R}^p} \ \frac{1}{2} \|\bgamma - \bbeta\|_2^2 + h ( \bgamma ). \]

The function $h$ is~\emph{essentially smooth} if $\inter(\dom(h))$ is nonempty, $h$ is differentiable on $\inter(\dom(h))$, and $\lim_{t \rightarrow + \infty} \|\nabla h \left( \bbeta^t \right) \|_2 = + \infty $ for any sequence $\{ \bbeta^t \} \subset \inter(\dom(h))$ converging to a boundary point of $\inter ( \dom (h) )$. 
The function $h$ is~\emph{locally smooth} if its gradient $\nabla h$ exists and is locally Lipschitz continuous on $\inter ( \dom (h) )$, that is, for every compact set $\mathcal{B} \subset \inter ( \dom (h) )$, there exists a constant $L_{\mathcal{B}} > 0$ such that $\| \nabla h( \bbeta ) - \nabla h \left( \bgamma \right) \|_2 \leq L_{\mathcal{B}} \| \bbeta - \bgamma \|_2$,
for all $\bbeta, \bgamma \in \mathcal{B}$.

The function $h$ is~\emph{piecewise linear-quadratic} (PLQ) if $\dom (h)$ can be represented as the union of finitely many polyhedral sets, relative to each of which $h$ is a quadratic function,~\textit{i.e.}, for each polyhedral set $\mathcal{P}$ in the union, there exist a symmetric matrix $\bm{Q} \in \mathbb{R}^{p \times p}$, a vector $\bb \in \mathbb{R}^p$, and a scalar $c \in \mathbb{R}$ such that $h(\bbeta) = \frac{1}{2} \bbeta^\top \bm{Q} \bbeta + \bb^\top \bbeta + c$ for all $\bbeta \in \mathcal{P}$.
The~\emph{subdifferential} of a proper convex function $h$ at a point $\bbeta \in \mathbb{R}^p$, denoted by $\partial h( \bbeta )$, is defined as the set of $\bs \in \mathbb{R}^p$ such that $ h \left( \bgamma \right) \geq h( \bbeta ) + \bs^\top \left( \bgamma - \bbeta \right)$ for all $\bgamma \in \mathbb{R}^p$.

A proper, closed, and convex function $h : \mathbb{R}^p \to \overline{\mathbb{R}}$ is said to be~\emph{firmly convex relative to a vector $\bnu \in \mathbb{R}^p$} if the tilted function $h_{\bnu} \left( \bbeta \right) := h( \bbeta ) - \bnu^\top \bbeta$ satisfies the~\emph{quadratic growth condition}.
That is, for any compact set $\mathcal{B} \subset \mathbb{R}^p$, there exists a constant $\kappa > 0$ such that
\begin{align*}
    h_{\bnu} \left( \bbeta \right) \geq \left( \inf h_{\bnu} \right) + \frac{\kappa}{2} \dist^2 \big( \bbeta, ( \partial h_{\bnu} )^{-1} ( \bm{0} )  \big) \qquad \forall \bbeta \in \mathcal{B}.
\end{align*}
The function $h$ is said to be~\emph{firmly convex} if $h$ is firmly convex relative to any $\bnu \in \mathbb{R}^p$ for which $\left( \partial h_{\bnu} \right)^{-1} \left( \bm{0} \right)$ is non-empty.

\subsection{Sparse GLMs, Branch-and-Bound, and Perspective Relaxations}
\label{sec:sparse_glm_problem_setup}
We aim to solve the following sparse GLM optimization problem to global optimality:
\begin{align} \label{obj:original_sparse_problem}
    \min\limits_{\bbeta \in \R^p} \left\{ f(\bX \bbeta, \by) + \lambda_2 \lVert \bbeta \rVert_2^2  \ : \ \| \bbeta \|_\infty \leq M, ~ \lVert \bbeta \rVert_0 \leq k \right\},
\end{align}
where $\bX \in \R^{n \times p}$ (with $\bX \neq \bm{0}$) and $\by \in \R^n$ denote the feature matrix and label vector, respectively.
Here, $M > 0$ is a user-specified box-constraint parameter, $f : \mathbb{R}^n \times \mathbb{R}^n \to \mathbb{R}$ is convex and differentiable, $k \in \mathbb{N}$ limits the number of nonzero coefficients, and $\lambda_2 > 0$ is a regularization parameter.

We use auxiliary binary variables $\bz \in \{0, 1\}^p$ and a big--M reformulation to arrive at the following mixed-integer nonlinear programming (MINLP) problem:
\begin{align}
    \label{obj:mip_formulation}
    \min_{ \left( \bbeta, \bz \right) \in \R^{2p} } \left\{ f(\bX \bbeta, \by) + \lambda_2 \lVert \bbeta \rVert_2^2 \ : \ \bz \in \{0,1\}^p, \ \bm{1}^\top \bz \leq k, \ |\beta_j| \leq M z_j, \ \forall j \in [p] \right\}.
\end{align}
At a high level, mixed-integer MIP solvers address such problems using branch-and-bound (BnB), recursively partitioning the search space by fixing some $z_j$ variables to $0/1$.

At each node $\mathcal N$ of the BnB tree, the perspective relaxation can be used to derive a tighter lower bound by solving the following convex relaxation:
\begin{align}
    \label{obj:perspective_relaxation}
    \min_{ \left( \bbeta, \bz \right) \in \mathcal{D}} \ f(\bX \bbeta, \by) + \lambda_2 \sum_{j \in [p]} \beta_j^2 / z_j,
\end{align}
where the domain $\mathcal{D}$ at node $\mathcal N$ is defined as
\begin{align}
    \label{eq:beta_z_domain}
    \mathcal{D} := \left\{
        \left( \bbeta, \bz \right) \in \mathbb{R}^{2p} \ \middle| \
        \begin{array}{l}
            \bm{1}^\top \bz \leq k, \quad \left| \beta_j \right| \leq M z_j, ~ \forall j \in [p], \\
            z_j = 0, ~ \forall j \in \mathcal{J}_0, \quad z_j = 1, ~ \forall j \in \mathcal{J}_1, \quad z_j \in [0, 1], ~ \forall j \in \mathcal{J}_f, \\
        \end{array}
    \right\}.
\end{align}
Here, the binary variables $\bz$ are relaxed to be continuous in $[0, 1]$, and the term $\beta_j^2/z_j$ is defined as $0$ if $\beta_j=z_j=0$ and $+\infty$ if $\beta_j \neq 0$ and $z_j=0$.
At any node $\mathcal N$ in the BnB tree, the sets $\mathcal{J}_0$, $\mathcal{J}_1$, and $\mathcal{J}_f$ form a partition of $[p]$ and correspond to the indices with $z_j=0$, $z_j=1$, and $z_j \in [0,1]$, respectively.
If $\mathcal{J}_0 = \mathcal{J}_1 = \emptyset$, then problem~\eqref{obj:perspective_relaxation} corresponds to the root node of the BnB tree.

\subsection{Composite Reformulation and Fenchel Duality}
\label{sec:problem_formulation}

We start by reformulating the perspective relaxation~\eqref{obj:perspective_relaxation} as an \emph{unconstrained} optimization problem. 
To this end, we define the implicit regularization function $g_{\mathcal{N}} : \mathbb{R}^p \rightarrow \overline{\mathbb{R}}$ as
\begin{align}
    \label{eq:function_g_definition}
    g_{\mathcal{N}}(\bbeta) := \inf \limits_{\bz \in \mathbb{R}^p} \left\{ \frac{1}{2} \sum_{j \in [p]} \beta_j^2 / z_j \ : \ (\bbeta, \bz) \in \mathcal{D} \right\},
\end{align}
where we adopt the convention that the infimum over an empty set is $+\infty$.
The function $g_{\mathcal{N}}$ captures the perspective term $\beta_j^2/z_j$ along with the cardinality and branching constraints on $\bz$ at node $\mathcal N$. 
We can therefore interpret $g_{\mathcal{N}}$ as a regularizer and rewrite~\eqref{obj:perspective_relaxation} in the following composite form:
\begin{align}
    \label{obj:original_sparse_problem_convex_composite_reformulation}
    \inf_{\bbeta \in \mathbb{R}^p} \big\{ \Phi\left(\bbeta\right):= F\left(\bX \bbeta\right) + G\left(\bbeta\right) \big\},
\end{align}
where $F\left(\bX \bbeta\right) := f\left(\bX \bbeta, \by\right)$ is the loss function and $G\left(\bbeta\right) := 2 \lambda_2 \, g_{\mathcal{N}}\left(\bbeta\right)$ is the regularization function.
Using the standard Fenchel duality framework, the dual problem of~\eqref{obj:original_sparse_problem_convex_composite_reformulation} is given by
\begin{align}
    \label{obj:original_sparse_problem_convex_composite_reformulation_fenchel_dual}
    \sup_{\bzeta \in \mathbb{R}^n} \left\{ \Psi\left(\bzeta\right):= -F^*\left(-\bzeta\right) - G^*\left(\bX^\top \bzeta\right) \right\},
\end{align}
where $F^*$ and $G^*$ are the Fenchel conjugates of $F$ and $G$, respectively.

By weak duality, $\Phi ( \bbeta ) \geq \Psi ( \bzeta )$ for any feasible pair $(\bbeta, \bzeta)$, meaning $\Psi( \bzeta )$ provides a valid lower bound for problem~\eqref{obj:original_sparse_problem}.
However, solving this problem efficiently is non-trivial. Directly applying standard first-order methods often yields slow sublinear convergence rates.
To address this, the next section introduces a \emph{general} primal--dual algorithmic framework.
By exploiting specific geometric properties of $F$ and $G$, we will show that this framework can achieve a linear convergence rate, enabling fast computation of tight lower bounds for the sparse GLM problems.

\section{A Linearly Convergent Algorithmic Framework}
\label{sec:methodology}
In this section, we develop a \emph{unified} analysis for composite optimization problems that satisfy specific geometric regularity conditions.
We consider the class of composite optimization problems defined by the primal and dual objective functions
\begin{align*}
\Phi \left( \bbeta \right) = F \left( \bX \bbeta \right) + G \left( \bbeta \right) \quad \text{ and }  \quad \Psi \left( \bzeta \right) = - F^* \left( - \bzeta \right) - G^* \left( \bX^\top \bzeta \right).
\end{align*}
We begin by stating the \emph{general regularity assumptions} required for our analysis.

\begin{assumption}  
    \label{assumption:F:G}
    The function $F : \mathbb{R}^n \to \mathbb{R}$ is convex and locally smooth. The function $G : \mathbb{R}^p \to \overline{\mathbb{R}}$ is proper, closed, and convex. Moreover, there are no lines along which $F$ or $G$ are finite and affine.
\end{assumption}

Since $F$ is assumed to be locally smooth, its domain necessarily has a nonempty interior. 
The assumption that there are no lines along which $F$ or $G$ are finite and affine further implies that the interiors of the domains of the conjugate functions $F^*$ and $G^*$ are nonempty \citep[Corollary~13.4.2]{rockafellar1970convex}.
We next define \( \Phi^\star := \inf_{\bbeta \in \mathbb R^p} \Phi(\bbeta) \) and \( \Psi^\star := \sup_{\bzeta \in \mathbb R^n} \Psi(\bzeta) \) as the optimal objective values of the primal and dual problems, respectively.
To support the primal--dual analysis, we impose the following Slater-type assumption.

\begin{assumption}
    \label{assumption:standard_fenchel_duality_assumption}
    The sets $\ri ( \dom (\Phi) )$ and $\inter ( \dom (\Psi) )$ are nonempty. Specifically, there exist Slater points $\bbeta_s \in \mathbb{R}^p$ and $\bzeta_s \in \mathbb{R}^n$ such that 
    $\bbeta_s \in \ri(\dom(G))$, $\bX \bbeta_s \in \inter(\dom(F))$, 
    $-\bzeta_s \in \inter(\dom(F^*))$, and $\bX^\top \bzeta_s \in \inter(\dom(G^*))$.
\end{assumption}

Under Assumption~\ref{assumption:standard_fenchel_duality_assumption}, the generalized Fenchel's Duality Theorem~\citep[Corollary 31.2.1]{rockafellar1970convex} ensures that both primal and dual problems are solvable and the strong duality holds. That~is, $$\Phi(\bbeta^\star) = \Phi^\star = \Psi^\star = \Psi(\bzeta^\star),$$
where $\bbeta^\star$ denotes an optimal primal solution and $\bzeta^\star$ denotes an optimal dual solution. 
We denote the set of all optimal primal solutions by $\mathcal{B}^\star$.
Furthermore, the KKT conditions~\citep[Theorem 31.3]{rockafellar1970convex} characterize any optimal primal--dual pair $(\bbeta^\star, \bzeta^\star)$ as
\begin{equation}
\label{eq:kuhn_tucker_conditions}
    \bX^\top \bzeta^{\star} \in \partial G(\bbeta^{\star}), \quad  \bX \bbeta^{\star} \in \partial F^*(-\bzeta^{\star}), \quad \bbeta^{\star} \in \partial G^*(\bX^\top \bzeta^{\star}), \quad  -\bzeta^{\star} = \nabla F(\bX \bbeta^{\star}).
\end{equation}
A key consequence of Assumption~\ref{assumption:F:G} is that $F$ is essentially smooth, which implies that its conjugate $F^*$ is essentially strictly convex \citep[Theorem 26.3]{rockafellar1970convex}, which ensures that the dual solution $\bzeta^\star$ is \emph{unique}.
We summarize these existence and uniqueness results in the following lemma.

\begin{lemma}
\label{lemma:basic}
    Under Assumptions~\ref{assumption:F:G} and \ref{assumption:standard_fenchel_duality_assumption}, the primal and dual problems are solvable and strong duality holds, that is, $\Phi^\star = \Psi^\star$. Moreover, the dual problem admits a unique optimal solution.
\end{lemma}

While Lemma~\ref{lemma:basic} guarantees the existence and uniqueness of solutions, these properties alone are insufficient to achieve linear convergence rates. To achieve faster convergence, the objective functions must exhibit specific curvature properties around the optimal set. In the following subsection, we formalize these requirements by introducing stronger geometric regularity conditions.

\subsection{Structural Geometry}

We now introduce the specific structural assumptions required to establish our main error bounds. These conditions strengthen standard convexity by imposing local curvature requirements compatible with the perspective relaxation structure.

\begin{assumption}
\label{assumption:quadratic}
    The following conditions hold.
    \begin{itemize}
        \item [(i)] The set of optimal primal solutions $\calB^\star$ is compact and satisfies $\bX \calB^\star \subset \inter(\dom (F))$.
        \item [(ii)] The functions $F$ and $G$ are firmly convex relative to $-\bzeta^\star$ and $\bX^\top \bzeta^\star$, respectively.
        \item [(iii)]
        The subdifferentials $\partial F^* ( - \bzeta^\star )$ and $\partial G^* ( \bX^\top \bzeta^\star )$ are polyhedral sets.
        \item [(iv)] The optimal dual solution satisfies $-\bzeta^\star \in \inter(\dom (F^*))$ and $\bX^\top \bzeta^\star \in \inter(\dom (G^*))$.
    \end{itemize}
\end{assumption}

Assumptions~\ref{assumption:quadratic}\,(i)--(iii) constitute the \emph{quadratic growth condition} commonly used in the primal analysis of first-order methods~\citep{zhou2017unified, drusvyatskiy2018error}. 
In contrast, Assumption~\ref{assumption:quadratic}\,(iv) targets the structural properties of the dual objective, a perspective that has received less attention in the existing literature. 
Under these combined assumptions, we establish the following fundamental geometric bounds.

\begin{theorem}
\label{theorem:quadratic:conditions}
    Under Assumptions~\ref{assumption:F:G}--\ref{assumption:quadratic}, let $\bzeta^\star$ denote the unique dual optimal solution.
    For any compact sets $\mathcal{B} \subset \mathbb{R}^p$ and $\mathcal{Z} \subset \mathbb{R}^n$ satisfying $\bzeta^\star \in \mathcal{Z}$, $-\mathcal{Z} \subset \inter(\dom(F^*))$ and $\bX^\top \mathcal{Z} \subset \inter(\dom(G^*))$, there exist constants $\alpha_{\mathcal{B}}, \kappa_{\mathcal{Z}} > 0$ such that
    \begin{align*}
        \Phi(\bbeta) \geq \Phi^\star + \frac{\alpha_{\mathcal{B}}}{2} \dist^2 (\bbeta, \mathcal{B}^\star), \, \forall \bbeta \in \mathcal{B}, \quad \text{and} \quad
        -\Psi(\bzeta) \leq -\Psi^\star + \frac{\kappa_{\mathcal{Z}}}{2} \|\bzeta - \bzeta^\star\|^2, \, \forall \bzeta \in \mathcal{Z}.
    \end{align*}
\end{theorem}

\begin{figure}[!htb]
    \vspace{0em}
    \centering
    \begin{tikzpicture}[scale=1.5, >=Stealth]

        \draw[thick] (-1.2, 1.4) parabola bend (-0.4,0) (0.4, 1.4) node[left] {$\Phi(\bbeta)$};

        \draw[dashed] (-1.6, 1.4) parabola bend (-0.4,0) (0.8, 1.4);

        \fill (-0.4, 0) circle (1.5pt) node[above, font=\small] {$\Phi^\star$};

        \draw[thick] (-0.8, -1.4) parabola bend (0.4,0) (1.6, -1.4) node[right] {$\Psi(\bzeta)$};

        \draw[dashed] (-0.3, -1.4) parabola bend (0.4,0) (1.1, -1.4);

        \fill (0.4, 0) circle (1.5pt) node[below, font=\small] {$\Psi^\star$};


        \node[right, align=left, font=\small] (growth) at (1.0, 0.8) {$\Phi^\star + \frac{\alpha_{\mathcal{B}}}{2} \dist^2 (\bbeta, \mathcal{B}^\star)$};
        \draw[->, shorten >=2pt] (growth.west) -- (0.45, 0.8);

        \node[below, align=center, font=\small] (decay) at (2.5, -0.56) {$\Psi^\star - \frac{\kappa_{\mathcal{Z}}}{2} \lVert \bzeta - \bzeta^\star \rVert^2$};
        \draw[->, shorten >=2pt] (decay.west) -- (0.9, -0.8);

    \end{tikzpicture}
    \caption{Illustration of quadratic growth for the primal objective $\Phi(\bbeta)$ and quadratic decay for the dual objective $\Psi(\bzeta)$.}
    \label{fig:quadratic_growth_decay}
    \vspace{0em}
\end{figure}
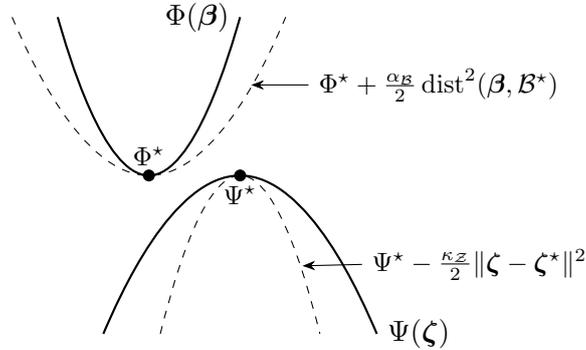

See Figure~\ref{fig:quadratic_growth_decay} for an illustration of the primal quadratic growth and dual quadratic decay conditions established in Theorem~\ref{theorem:quadratic:conditions}.
Theorem~\ref{theorem:quadratic:conditions} reveals a striking symmetry.
While the primal objective function $\Phi(\bbeta)$ \emph{grows at least} quadratically away from the optimal solution set $\mathcal{B}^\star$, the dual objective function $\Psi(\bzeta)$ \emph{decays at most} quadratically away from the unique dual solution $\bzeta^\star$. 
Consequently, we refer to the property guaranteed by Assumption~\ref{assumption:quadratic}\,(iv) as the \emph{quadratic decay condition}, characterizing it as the dual counterpart to the primal quadratic growth condition presented in Assumptions~\ref{assumption:quadratic}\,(i)--(iii).

The proof of Theorem~\ref{theorem:quadratic:conditions} relies on two technical lemmas regarding firmly convex functions. The first establishes a subdifferential error bound.

\begin{lemma}
    \label{lemma:error_bound_for_fenchel_conjugate_of_firmly_convex_function}
    Let $h$ be a proper, closed, convex function such that there are no lines along which $h$ is finite and affine.
    If $h$ is firmly convex relative to some $\bnu \in \inter(\dom (h^*))$, then for any compact set $\mathcal U \subset \inter(\dom (h^*))$ with $\bnu \in \mathcal U$, there exists a constant $\alpha_{\mathcal U} > 0$ such that
    \begin{align*}
        \dist \left( \partial h^*\left( \bmu \right), \partial h^*\left( \bnu \right)   \right) \leq \frac{2}{\alpha_{\mathcal{U}}} \|\bmu - \bnu \|_2, \quad ~\forall \bmu \in \mathcal{U}.
    \end{align*}
\end{lemma}

The proof of Lemma~\ref{lemma:error_bound_for_fenchel_conjugate_of_firmly_convex_function} is provided in Appendix~\ref{ec_proof:error_bound_for_fenchel_conjugate_of_firmly_convex_function}.
Lemma~\ref{lemma:error_bound_for_fenchel_conjugate_of_firmly_convex_function} can be viewed as the dual counterpart to the subdifferential error bound in \citep[Theorem~3.3]{drusvyatskiy2018error}. 
While \citet{drusvyatskiy2018error} analyze the primal subdifferential $\partial h$, our result characterizes the behavior of the conjugate subdifferential $\partial h^*$. 
Notably, the inequality direction is reversed, reflecting the sign flip inherent in the Fenchel dual formulation.
Building on this error bound, we derive a quadratic upper bound for the conjugate function value.

\begin{lemma}
    \label{lemma:quadratic_upper_bound_for_fenchel_conjugate_of_firmly_convex_function}
    Let $h$ be a proper, closed, convex function such that there are no lines along which $h$ is finite and affine.
    If $h$ is firmly convex relative to some $\bnu \in \inter(\dom (h^*))$, then for any compact set $\mathcal U \subset \inter(\dom (h^*))$ with $\bnu \in \mathcal U$, there exists a constant $\alpha_{\mathcal U} > 0$ and a vector $\bs \in \partial h^*(\bnu)$ such that
    \begin{align*}
        h^*(\bmu) \leq h^*(\bnu) + \bs^\top \left( \bmu - \bnu\right) + \frac{\alpha_{\mathcal{U}}}{2} \|\bmu - \bnu\|^2, \quad ~\forall \bmu \in \mathcal{U}.
    \end{align*}
\end{lemma}

The proof of Lemma~\ref{lemma:quadratic_upper_bound_for_fenchel_conjugate_of_firmly_convex_function} is provided in Appendix~\ref{ec_proof:quadratic_upper_bound_for_fenchel_conjugate_of_firmly_convex_function}.
We can now prove Theorem~\ref{theorem:quadratic:conditions} by leveraging Lemma~\ref{lemma:quadratic_upper_bound_for_fenchel_conjugate_of_firmly_convex_function}.
The proof of Theorem~\ref{theorem:quadratic:conditions} is provided in Appendix~\ref{ec_proof:quadratic_conditions}.
The primal quadratic growth and dual quadratic decay established above constitute the fundamental geometric properties of our problem structure.

\subsection{Primal--Dual Error Bounds}

While Theorem~\ref{theorem:quadratic:conditions} characterizes the primal and dual landscapes \emph{independently}, we now establish the rigorous connection between them.
To facilitate this analysis, we first introduce the following regularity condition on the primal loss function.

\begin{assumption}
\label{assumption:loss_function_is_essentially_strictly_convex}
    The loss function $F$ is strictly convex.
\end{assumption}


With this assumption in place, we can relate the primal optimality gap $\Phi(\bbeta) - \Phi^\star$ to two critical dual quantities: the dual distance to optimality $\|\bzeta - \bzeta^\star\|_2^2$ and the dual optimality gap $\Psi^\star - \Psi(\bzeta)$. 
The link between these quantities relies on a specific construction of the dual variable. 
Inspired by the KKT condition $-\bzeta^\star = \nabla F(\bX \bbeta^\star)$ in~\eqref{eq:kuhn_tucker_conditions}, we define the \emph{induced dual vector} $\bzeta$ associated with an arbitrary primal feasible solution $\bbeta$ via the mapping
\[
    \bzeta := - \nabla F(\bX \bbeta).
\]
In the following theorem, we prove that under this mapping, the primal optimality gap strictly controls the dual error, effectively allowing the primal convergence to drive the dual convergence.

\begin{theorem}
    \label{theorem:primal_dual_relation_dual_sequence_and_function_value}
    Under Assumptions~\ref{assumption:F:G}--\ref{assumption:loss_function_is_essentially_strictly_convex}, let $\mathcal{B}^\star$ denote the set of primal optimal solutions and
    $\bzeta^\star$ denote the unique dual optimal solution.
    Let $\mathcal{B} \subseteq \dom (G)$ be a compact set such that $\mathcal{B}^\star \subseteq \mathcal{B}$. 
    Define the corresponding dual image set $\mathcal{Z} := \{ - \nabla F(\bX \bbeta) \mid \bbeta \in \mathcal{B} \}$, and assume that the interior conditions $-\mathcal{Z} \subset \inter(\dom (F^*))$ and $\bX^\top \mathcal{Z} \subset \inter(\dom (G^*))$ hold. 
    Then, there exist constants $\sigma_{\mathcal{B}}, \kappa_{\mathcal{B}} > 0$ such that for any primal vector $\bbeta \in \mathcal{B}$ and its induced dual vector $\bzeta := - \nabla F(\bX \bbeta)$, we have
    \begin{align*}
        \frac{\sigma_{\mathcal{B}}}{2 } \|\bzeta - \bzeta^\star\|_2^2 &\leq \Phi(\bbeta) - \Phi^\star  \quad 
        \text{and} \quad 
        \Psi(\bzeta^\star) - \Psi(\bzeta) \leq \frac{\kappa_{\mathcal{B}}}{\sigma_{\mathcal{B}}} \big( \Phi(\bbeta) - \Phi^\star \big). 
    \end{align*}
\end{theorem}

The proof of Theorem~\ref{theorem:primal_dual_relation_dual_sequence_and_function_value} is provided in Appendix~\ref{ec_proof:primal_dual_relation_dual_sequence_and_function_value}.
To provide geometric intuition for the first inequality in Theorem~\ref{theorem:primal_dual_relation_dual_sequence_and_function_value}, we draw a parallel to standard results in smooth convex optimization.
Recall that if a function $h$ is globally $L$-smooth, its Fenchel conjugate $h^*$ is $(1/L)$-strongly convex~\citep[Theorem 5.10]{beck2017first}. 
Moreover, for any $\bz \in \mathbb{R}^p$ and a minimizer $\bz^\star = \argmin h(\bz)$, the smoothness of $h$ implies a lower bound on the optimality gap~\citep[Theorem 5.8]{beck2017first}:
\begin{align*}
    \frac{1}{2L} \| \nabla h(\bz) \|_2^2 \leq h(\bz) - h(\bz^\star).
\end{align*}
If we interpret $\nabla h(\bz)$ as a dual variable, this inequality is analogous to the first bound in Theorem~\ref{theorem:primal_dual_relation_dual_sequence_and_function_value}.
However, our result extends this relationship to a more complex setting where the primal objective $F+G$ has a non-smooth composite structure, and the loss function $F$ is only \emph{locally} smooth.

Moreover, Theorem~\ref{theorem:primal_dual_relation_dual_sequence_and_function_value} has a crucial algorithmic implication: the convergence of the primal objective $\Phi(\bbeta^t)$ strictly controls the convergence of both the dual sequence $\bzeta^t$ and the dual objective $\Psi(\bzeta^t)$. 
This relationship holds independently of the specific algorithm used to generate the primal sequence. 
Consequently, if an algorithm ensures that $\Phi(\bbeta^t)$ converges to $\Phi^\star$ at a specific rate, the dual sequence and dual objective values inherit the same rate.






\subsection{Linear Convergence via Gap-Based Restarts}
\label{subsec:restart_scheme}

To derive a unified analysis that applies to a wide range of first-order methods, we first formalize the properties of a generic solver. Rather than tailoring our analysis to the specific update rules of a single algorithm, we assume access to a solver that generates iterates within a compact region and provides a standard convergence guarantee.

\begin{assumption}
    \label{assumption:sublinear_algorithm}
    Let $\mathcal{B} \subseteq \dom G$ be a compact set such that $\mathcal{B}^\star \subseteq \mathcal{B}$.
    There exists a proximal first-order algorithm that generates a feasible primal sequence $\{\bbeta^t\}_{t \ge 0} \subset \mathcal{B}$ and an associated dual sequence $\bzeta^t := -\nabla F(\bX \bbeta^t)$.
    We assume that the algorithm relies on the proximal operator of $G$.
    Furthermore, the algorithm guarantees a convergence rate relative to the initialization $\bbeta^0$ characterized by a rate function $r: \mathbb{N} \to \mathbb{R}_+$, that is,
    \begin{equation}
        \label{eq:algorithm_rate_assumption}
        \Phi(\bbeta^t) - \Phi^\star \leq r(t) \cdot \dist^2(\bbeta^0, \mathcal{B}^\star) \quad \forall t \ge 1.
    \end{equation}
\end{assumption}

Assumption~\ref{assumption:sublinear_algorithm} is satisfied by a broad class of first-order methods, including standard PGD \citep[Section 10.2]{beck2017first}, accelerated variants such as FISTA \citep[Section 10.7]{beck2017first}, and parameter-free methods including AdProxGD \citep{malitsky2024adaptive} and AC-FGM \citep{li2023simple}.
We first observe that any algorithm satisfying this assumption automatically yields a convergence guarantee for the duality gap. 
By combining the algorithmic rate in~\eqref{eq:algorithm_rate_assumption} with the structural bounds established in Theorem~\ref{theorem:primal_dual_relation_dual_sequence_and_function_value}, we derive the following duality gap.

\begin{proposition}
    \label{prop:duality_gap_sublinear_bound}
    Under Assumptions~\ref{assumption:F:G}--\ref{assumption:sublinear_algorithm}, let $\alpha_{\mathcal{B}}, \sigma_{\mathcal{B}}, \kappa_{\mathcal{B}} > 0$ be the geometric constants derived in Theorem~\ref{theorem:primal_dual_relation_dual_sequence_and_function_value}.
    Then, the primal--dual iterates $\{(\bbeta^t, \bzeta^t)\}_{t \ge 1}$, satisfy the duality gap 
    \begin{align}
        \label{eq:duality_gap_contraction}
        \Phi(\bbeta^t) - \Psi(\bzeta^t) \leq \left( \frac{2 r(t) (1 + \kappa_{\mathcal{B}} / \sigma_{\mathcal{B}})}{\alpha_{\mathcal{B}}} \right) \big( \Phi(\bbeta^0) - \Psi(\bzeta^0) \big).
    \end{align}
\end{proposition}
The proof of Proposition~\ref{prop:duality_gap_sublinear_bound} is provided in Appendix~\ref{ec_proof:duality_gap_sublinear_bound}.
Proposition~\ref{prop:duality_gap_sublinear_bound} establishes the duality gap shares the same convergence rate as $r(t)$, which is sublinear for standard first-order methods.
A notable exception is PGD with a fixed step-size, which automatically achieves a linear primal convergence rate under the quadratic growth condition~\citep{zhou2017unified, drusvyatskiy2018error}.
This is consistent with Proposition~\ref{prop:duality_gap_sublinear_bound} because it immediately implies a linear convergence rate for PGD on the duality gap.

In contrast, accelerated methods like FISTA typically require restart schemes to prevent oscillations and unlock linear convergence for primal sequence~\citep{o2015adaptive, fercoq2019adaptive}.
Furthermore, recent adaptive methods~\citep{malitsky2024adaptive, li2023simple} generally yield only sublinear rates, and their theoretical guarantees for the dual objective under quadratic growth remain largely unexplored.
To ensure linear convergence (in the standard Q-linear sense)
across this broader class of algorithms, we employ a restart strategy governed explicitly by the duality gap.
The core idea is to run the algorithm only until the duality gap decreases by a fixed factor $\eta > 1$, and then restart the algorithm, treating the current iterate as the new initialization. 
Because the~duality gap is a strict upper bound on the error, this ensures the error contracts geometrically at each restart~epoch.

\begin{theorem}
    \label{theorem:restart_linear_convergence}
    Suppose the conditions of Proposition~\ref{prop:duality_gap_sublinear_bound} hold and $r\left( t \right)$ is monotonically decreasing.
    Let $(\bbeta^{\init}, \bzeta^{\init})$ be the initial primal--dual pair.
    For any target contraction factor $\eta > 1$, define the restart interval $t_{\max}$ as
    \begin{equation}
        \label{eq:tmax}
        t_{\max} := \left\lceil r^{-1} \left( \frac{ \alpha_{\mathcal{B}} }{ 2 \eta ( 1 + \kappa_{\mathcal{B}} / \sigma_{\mathcal{B}} ) } \right) \right\rceil,
    \end{equation}
    where $r^{-1}$ is the inverse of the rate function.
    If the algorithm is restarted every $t_{\max}$ iterations (resetting $t=0$ and setting $\bbeta^0$ to the current iterate), the sequence of duality gaps at the restart points converges linearly with
    \begin{align*}
        \Phi(\bbeta^{s \cdot t_{\max}}) - \Psi(\bzeta^{s \cdot t_{\max}}) \leq \eta^{-s} \left( \Phi(\bbeta^{\init}) - \Psi(\bzeta^{\init}) \right).
    \end{align*}
\end{theorem}

The proof of Theorem~\ref{theorem:restart_linear_convergence} is provided in Appendix~\ref{ec_proof:restart_linear_convergence}.
A crucial consequence of Theorem~\ref{theorem:restart_linear_convergence} is that the linear convergence of the duality gap necessitates the linear convergence of the primal and dual objective values, as well as their respective iterates.
To see this, note that by weak duality ($\Psi(\bzeta) \leq \Psi^\star = \Phi^\star$), the duality gap serves as a strict upper bound for the primal optimality gap
\begin{align*}
    \Phi \left(\bbeta^{s \cdot t _{ \text{max} }} \right) - \Phi^\star
    \leq \Phi \left( \bbeta^{s \cdot t _{ \text{max} }} \right) - \Psi \left( \bzeta^{s \cdot t _{ \text{max} }} \right)
    \leq \eta^{-s} \left( \Phi(\bbeta^0) - \Psi(\bzeta^0) \right).
\end{align*}
This establishes the linear convergence of primal objective values.
By Theorems~\ref{theorem:quadratic:conditions} and~\ref{theorem:primal_dual_relation_dual_sequence_and_function_value}, primal iterates $\{\bbeta^t\}$, dual iterates $\{\bzeta^t\}$, and dual objective values $\{\Psi(\bzeta^t)\}$ also converge linearly.

We conclude this section by quantifying the specific iteration complexities for PGD, FISTA, and AC-FGM. 
For AC-FGM, the standard sublinear rate depends on user-specified parameters $c \in (1, \infty)$ and $\rho \in (0, 1 - \sqrt{6}/3]$, where $c$ is employed for the line-search procedure during the initialization. 
Additionally, $L_0$ denotes an initial underestimate of the smoothness constant used by the algorithm.
To simplify notation in the comparison below, we omit the subscript $\mathcal{B}$ from the geometric constants $\kappa$ (dual decay), $L$ (primal smoothness, \textit{with respect to} $F(\bX \cdot)$), $\alpha$ (primal growth), and $\sigma$ (dual strong convexity).
Table~\ref{table:restart_complexity} summarizes the rate functions $r(t)$ and the resulting complexities after applying our duality-gap-based restart scheme.
Crucially, the restart scheme transforms the worst-case sublinear convergence of the accelerated methods (FISTA and AC-FGM) into a linear rate that scales with the \emph{square root} of the condition number ($\sqrt{\kappa L / \alpha \sigma}$). 
This represents a significant theoretical improvement over the standard PGD linear rate, which scales linearly with the condition~number.

\begin{table}[!htb]
    \caption{Standard rate functions $r(t)$ (when only the smoothness of $F$ is assumed) and the corresponding complexity to reach $\epsilon$-accuracy after applying the duality-gap-based restart scheme.}
    \centering
    \label{table:restart_complexity}
    \begin{tabular}{lcc}
        {Algorithm} & {Rate Function} $r(t)$ & {Complexity with Restart} \\
        \hline
        PGD & $\displaystyle \frac{L}{2t}$ & $\displaystyle \mathcal{O} \left( \frac{ \kappa L }{\alpha \sigma} \log\left( \frac{1}{\epsilon} \right) \right)$ \\
        \hline
        FISTA & $\displaystyle \frac{2L}{\left( t+1 \right) ^{ 2 }}$ & $\displaystyle \mathcal{O} \left( \sqrt{\frac{ \kappa L }{\alpha \sigma }} \log\left( \frac{1}{\epsilon} \right)\right)$ \\
        \hline
        AC-FGM & $\displaystyle \frac{15 c L}{2 t \left( t+1 \right) \rho \left( 1 - \rho \right) } $ & $\displaystyle \mathcal{O} \left( \sqrt{\frac{ \kappa L }{\alpha \sigma}} \log\left( \frac{1}{\epsilon} \right) + \log_{c} \left( \frac{L}{L_0}  \right) \right)$
    \end{tabular}
\end{table}

\subsection{Restart in Practice}
Theorem~\ref{theorem:restart_linear_convergence} provides a theoretical guarantee based on the geometric constants $\sigma_{\mathcal{B}}$, $\alpha_{\mathcal{B}}$, and $\kappa_{\mathcal{B}}$, as well as the specific rate function $r(t)$. 
In practice, these parameters are often unknown or difficult to estimate.
However, it is not required to know these constants to apply the restart scheme.
Because the duality gap is computable, we can monitor it at runtime and trigger a restart dynamically.

Specifically, the practical implementation involves checking the condition $\Phi(\bbeta^t) - \Psi(\bzeta^t) \leq \eta^{-1} (\Phi(\bbeta^0) - \Psi(\bzeta^0))$ at each iteration.
Theoretical analysis ensures that this condition will be met within at most $t_{\max}$ iterations, where $t_{\max}$ is defined as in~\eqref{eq:tmax}.
To achieve an $\epsilon$-accuracy for the primal--dual gap, the total number of restarts $s_{\text{max}}$ required is upper bounded by
\[
    s_{\text{max}} \leq \left\lceil \log_\eta \left( \frac{\Phi(\bbeta^{\init}) - \Psi(\bzeta^{\init})}{\epsilon} \right) \right\rceil,
\]
which is logarithmic in $1/\epsilon$.
Thus, the total number of iterations is upper bounded by
\[
    T_{\text{total}} \leq s_{\text{max}} \cdot t_{\text{max}} \leq \left\lceil \log_\eta \left( \frac{\Phi(\bbeta^{\text{init}}) - \Psi(\bzeta^{\text{init}})}{\epsilon}\right) \right\rceil \cdot \left\lceil r^{-1} \left( \frac{ \alpha_{\mathcal{B}} }{ 2 \eta ( 1 + \kappa_{\mathcal{B}} / \sigma_{\mathcal{B}} ) } \right) \right\rceil.
\]
The following theorem gives the optimal restart parameter $\eta$ under a polynomial rate model.
\begin{theorem}
    \label{theorem:optimal_restart_eta_choice}
    Suppose the convergence rate behaves like $r(t) \approx C / t^q$ for some $q > 0$ and $C > 0$.
    Then, the constant-optimal choice of the restart parameter is $\eta = e^q$.
\end{theorem}

The proof of Theorem~\ref{theorem:optimal_restart_eta_choice} is provided in Appendix~\ref{ec_proof:optimal_restart_eta_choice}.
In particular, $q=2$ (typical for accelerated methods, where $r(t)=\Theta(1/t^2)$) yields $\eta=e^2$.
However, note that this optimal choice is derived from a worst-case scenario and may not match the best empirical choice of $\eta$.
For example, it is possible to observe that $\eta=e^3$ performs better than $\eta=e^2$ in practice, for at least three reasons.
First, the bound $r(t)=\Theta(1/t^2)$ for accelerated methods is typically not tight: if the effective decrease within an epoch behaves like $r(t)\approx C/t^q$ with some $q>2$ on the problem instances of interest, then the same calculation above suggests a larger $\eta$ (namely $\eta=e^q$) is preferred.
Second, our derivation treats each restart as ``free'' beyond the iteration budget $t_{\max}$.
In practice, restarting an accelerated method discards its momentum and introduces a transient ``warm-up'' period before the accelerated regime becomes effective.
If we model this effect as an additional overhead of $t_0$ iterations per restart (or, more generally, a restart-dependent constant in the bound), then larger $\eta$ reduces the number of restarts and therefore pays the overhead fewer times.
Third, our practical scheme triggers restarts based on the \emph{observed} duality-gap decrease rather than the worst-case bound: if the gap often drops by a factor of $\eta^{-1}$ much earlier than $t_{\max}$, then increasing $\eta$ can reduce the number of restarts substantially while increasing the typical epoch length only mildly.
Taken together, these effects explain why $\eta=e^2$ is best viewed as the optimal choice under an idealized polynomial-rate model, whereas larger values such as $\eta=e^3$ can be preferable empirically.

\section{Efficient Implementation for Sparse GLMs}

The algorithmic framework in Section~\ref{sec:methodology} guarantees linear convergence for generic composite problems,
provided they satisfy Assumptions~\ref{assumption:F:G}--\ref{assumption:loss_function_is_essentially_strictly_convex}. 
In this section, we apply this framework to the perspective relaxation of sparse GLMs at each node $\mathcal N$ of the BnB tree.

Recall from Section~\ref{sec:problem_formulation} that this problem takes the composite form with $F(\bX \bbeta) := f(\bX \bbeta, \by)$ and $G(\bbeta) := 2 \lambda_2 g_{\mathcal{N}}(\bbeta)$, where $g_{\mathcal{N}}$ is the implicit regularizer defined as
\begin{align*}
    g_{\mathcal{N}}(\bbeta) = \inf \limits_{\bz \in \mathbb{R}^p} \left\{ \frac{1}{2} \sum_{j \in [p]} \beta_j^2 / z_j \ : \ (\bbeta, \bz) \in \mathcal{D} \right\}.
\end{align*}
Here, the domain $\mathcal{D}$, defined in~\eqref{eq:beta_z_domain}, encodes the branching constraints at node $\mathcal{N}$ by specifying which binary variables are fixed to $0$ or $1$, and which remain free.

To fully leverage the theoretical guarantees and practical performance of our framework, we address three critical aspects in the following subsections.
First, we establish that $g_{\mathcal{N}}$ possesses the necessary geometric properties.
Second, we demonstrate that the proximal operator of $g_{\mathcal{N}}$ can be evaluated efficiently, ensuring that the per-iteration computational cost remains low.
Finally, we verify that our primal--dual framework can be applied by showing that standard GLM loss functions naturally satisfy the requisite smoothness and strict convexity assumptions.

\subsection{Geometric Properties of the Implicit Regularizer}

We begin by establishing the fundamental geometric properties of the implicit regularizer $g_{\mathcal{N}}$.
First, we verify that $g_{\mathcal{N}}$ is a valid convex regularizer with a bounded domain, which is a prerequisite for the compactness of the optimal solution set.

\begin{lemma}
\label{lemma:g_is_proper_closed_and_convex}
    The implicit function $g_{\mathcal{N}}$ in~\eqref{eq:function_g_definition} is proper, closed and convex, with a compact domain.
\end{lemma}

The proof of Lemma~\ref{lemma:g_is_proper_closed_and_convex} is provided in Appendix~\ref{ec_proof:g_is_proper_closed_and_convex}.
The next lemma gives a closed-form expression for $g^*$, where the function $\TopSum_{k}(\cdot)$ denotes the sum of the top $k$ largest
elements, and $\bm{H}_M \left( \balpha \right) := \left( H_M \left( \alpha_1 \right), \ldots, H_M \left( \alpha_p \right)    \right) $ is the Huber loss function applied to each coordinate with
\begin{align}
    \label{eq:huber_loss_in_bnb}
    & H_M(\alpha_j) := \begin{cases}
        \frac{1}{2} \alpha_j^2 & \text{if } \lvert \alpha_j \rvert \leq M \\
        M \lvert \alpha_j \rvert - \frac{1}{2} M^2 & \text{if } \lvert \alpha_j \rvert > M.
    \end{cases}
\end{align} 
\begin{lemma}
    \label{lemma:function_g_conjugate_in_BnB}
    The Fenchel conjugate of $g_{\mathcal{N}}$ is given by 
    \begin{align}
        \label{eq:fenchel_conjugate_of_function_g_in_bnb}
        g_{\mathcal{N}}^* \left( \balpha \right) = 
            \sum_{j \in \mathcal{J}_1} H_M(\alpha_j) \, + \textstyle
            \TopSum_{k - | \mathcal J_1|} \big( \bm{H}_M \big( \balpha_{\calJ_f} \big) \big).
    \end{align}
\end{lemma}

The proof of Lemma~\ref{lemma:function_g_conjugate_in_BnB} is provided in Appendix~\ref{ec_proof:function_g_conjugate_in_BnB}.
The second term in~\eqref{eq:fenchel_conjugate_of_function_g_in_bnb} reveals that when we shift from the \emph{primal} space to the \emph{dual} space, the variational formulation of $g_{\mathcal{N}}$ is transformed into an explicit expression.
The closed-form expression in Lemma~\ref{lemma:function_g_conjugate_in_BnB} reveals that $g_{\mathcal{N}}^*$ is constructed entirely from operations, namely finite sums, sorting, and quadratic splines, that preserve the PLQ structure.
This observation is formalized in the following lemma.

\begin{lemma}
    \label{lemma:g_and_g_Fenchel_are_PLQ_and_firmly_convex}
    The functions $g_{\mathcal{N}}$ and $g_{\mathcal{N}}^*$ are both PLQ and firmly convex.
\end{lemma}

\begin{figure}[!htb]
    \vspace{0em}
    \centering
    \includegraphics[width=0.98\linewidth]{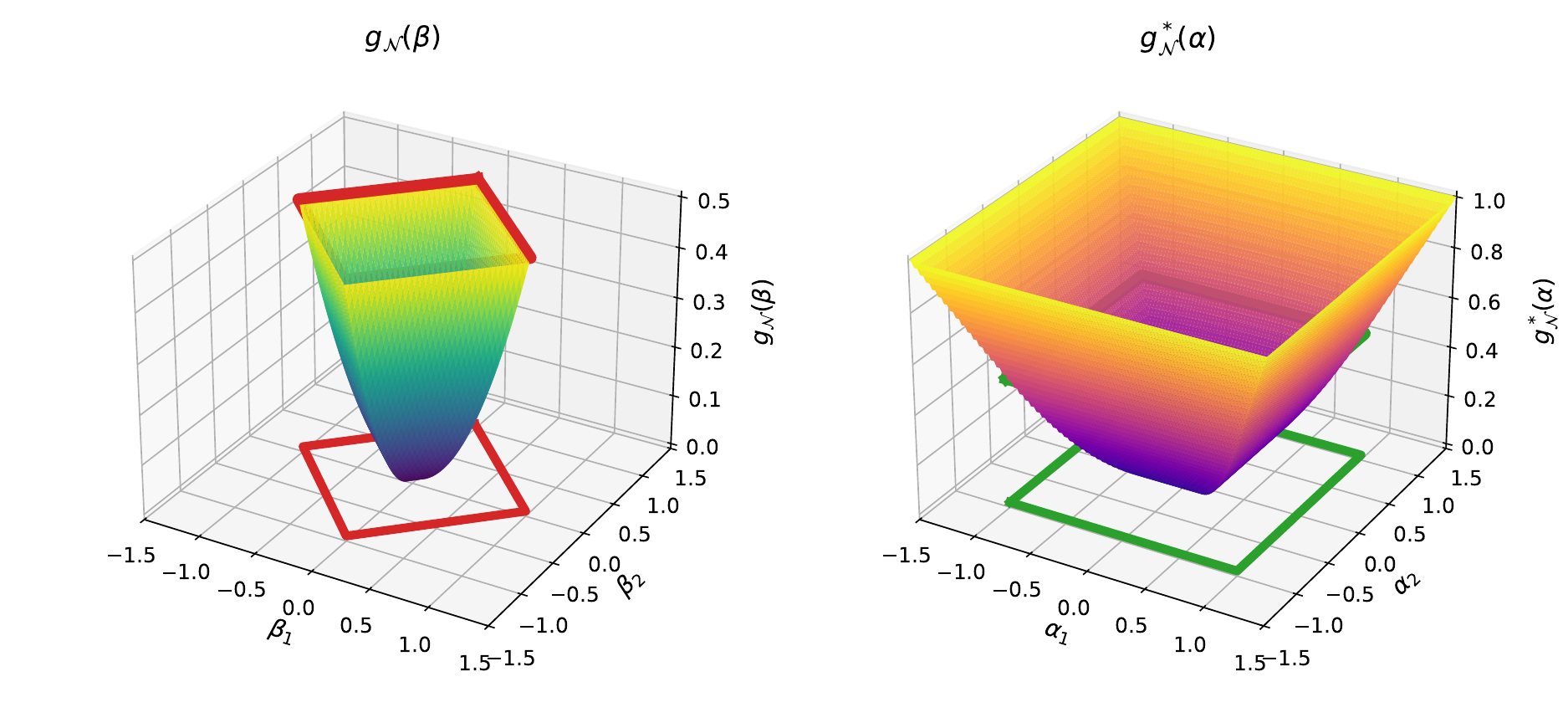}
    \caption{3D plots of $g_{\mathcal{N}}$ (left) and $g_{\mathcal{N}}^*$ (right) at the root node ($p=2$, $k=1$, $M=1$). Left: red diamonds mark the boundary $\lvert \beta_1 \rvert + \lvert \beta_2 \rvert = M$ of $\dom(g_{\mathcal{N}})$ (top at $g_{\mathcal{N}}(\bbeta)=\frac{1}{2}M^2$, bottom projection). Right: green squares mark the transition set $\max(\lvert \alpha_1 \rvert, \lvert \alpha_2 \rvert)=M$ (top at $g_{\mathcal{N}}^*(\balpha)=\frac{1}{2}M^2$, bottom projection); $g_{\mathcal{N}}^*$ grows quadratically inside and linearly outside.}
    \label{fig:gN_and_gNstar_plq_3d}
    \vspace{0em}
\end{figure}

The proof of Lemma~\ref{lemma:g_and_g_Fenchel_are_PLQ_and_firmly_convex} is provided in Appendix~\ref{ec_proof:g_and_g_Fenchel_are_PLQ_and_firmly_convex}.
This result confirms that the perspective regularizer satisfies the rigorous geometric requirements of our framework. 
See Figure~\ref{fig:gN_and_gNstar_plq_3d} for a visual illustration of the PLQ structures of both $g_{\mathcal{N}}$ and its conjugate $g_{\mathcal{N}}^*$.
We will leverage this property in Proposition~\ref{proposition:regularity} within Section~\ref{sec:conv:rate:GLM}.

\subsection{Efficient Evaluation of the Implicit Regularizer and its Proximal Operator}

Having established the geometric regularity of the implicit regularizer, we now address its computational tractability.
Our algorithmic framework requires two key operations: (i) evaluating the function value $g_{\mathcal{N}}(\bbeta)$ in order to monitor the duality gap for restarts and (ii) computing its proximal operator in the proximal gradient method.

We first address function evaluation. 
Note that directly solving the minimization problem in~\eqref{eq:function_g_definition} via generic Second-Order Cone Programming (SOCP) solvers is computationally expensive and, due to the iterative nature of interior-point methods, yields only approximate solutions.
In contrast, by exploiting the specific structure of the domain $\mathcal{D}$, we derive a specialized algorithm that~evaluates $g_{\mathcal{N}}(\bbeta)$ \emph{exactly} in a finite number of steps, with a worst-case computational complexity log-linear~in~$p$.
The overall algorithm is presented in Algorithm~\ref{alg:compute_g_value_root_node_algorithm}.

\begin{theorem}
    \label{theorem:compute_g_value_algorithm_correctness}
    For any $\bbeta \in \dom (g_{\mathcal{N}})$, let $\bar \bbeta := \bbeta_{\mathcal{J}_f}$, $\bar p := \lvert \mathcal{J}_f \rvert$, and $\bar k := k - \lvert \mathcal{J}_1 \rvert$.
    Then,
    \begin{align*}
        g_{\mathcal{N}}(\bbeta)
        &= \frac{1}{2} \sum_{j \in \mathcal{J}_1}\beta_j^2 + \frac{1}{2}\min_{\bomega} \left\{ \sum_{j=1}^{\bar k} \omega_j^2 \ \middle|\ 
        \begin{array}{l}
            \bomega \ge \bm{0}, \quad M \ge \omega_1 \ge \cdots \ge \omega_{\bar k} \ge 0, \\
            \omega_{\bar k+1} = \cdots = \omega_{\bar p} = 0, \quad \bomega \succeq_m |\bar{\bbeta}|
        \end{array}
        \right\}.
    \end{align*}
    Here, $\bomega \succeq_m |\bar{\bbeta}|$ denotes that $\bomega$ majorizes $|\bar{\bbeta}|$, \textit{i.e.},
    \[
        \sum_{j=1}^{\ell} \omega_j \ge \sum_{j=1}^{\ell} |\bar\beta_{(j)}| \quad \forall \ell = 1,\ldots,\bar p-1,
        \qquad \text{and} \qquad
        \sum_{j=1}^{\bar p} \omega_j = \sum_{j=1}^{\bar p} |\bar\beta_{(j)}|.
    \]
    Algorithm~\ref{alg:compute_g_value_root_node_algorithm} computes a minimizer $\bomega$ of the above problem and thus evaluates $g_{\mathcal{N}}(\bbeta)$ exactly with a computational complexity of $\mathcal O(p + \bar p \log \bar k + \bar k)$.
\end{theorem}

\begin{algorithm}[!tb]
    \DontPrintSemicolon
        \caption{Compute $g_{\mathcal{N}}(\bbeta)$ at node $\mathcal{N}$ with associated index sets $\mathcal{J}_0,  \mathcal{J}_1, \mathcal J_f$}
        \label{alg:compute_g_value_root_node_algorithm}
        \KwData{vector $\bbeta \in \dom (g_{\mathcal{N}})$, cardinality parameter $k \in [p]$}
        \KwResult{$g_{\mathcal{N}}(\bbeta)$}
        
        Set $\bar k \gets k - | \mathcal J_1 |, \ \bar p \leftarrow | \mathcal J_f |$ \;
        Set $\bar \bbeta \gets \bbeta_{\mathcal J_f},  \ \bomega \gets \bm{0}_{\bar k}, \ \theta \gets \sum_{j \in [\bar p]} \vert{\bar \beta_j} $ \;
        Sort $\bar \bbeta$ such that $\vert{\bar \beta_{(1)}} \geq \ldots \geq \vert{\bar \beta_{(\bar k)}} \geq \max\limits_{j \in [\bar k+1 : \bar p]} \{ \vert{\bar \beta_j} \}$\;
        \For{$j \gets 1$ \KwTo $\bar k$}{
            $\overline{\theta} \gets \theta / (\bar k - j + 1)$ \;
            \lIf{$\overline{\theta} \geq \vert{\bar \beta_{(j)}}$}
            {
                $\omega_{[j:\bar k]} \gets \overline{\theta}$ \textbf{break}
            }
            \lElse{
                $\omega_j \gets \vert{\bar \beta_{(j)}}$ and $\theta \gets \theta - \vert{\bar \beta_{(j)}}$
            }
        }
        \KwRet{$\frac{1}{2} \sum_{j \in \mathcal{J}_1} \beta_j^2 +  \frac{1}{2}\sum_{j \in [\bar k]} \omega_j^2$}\;
\end{algorithm}

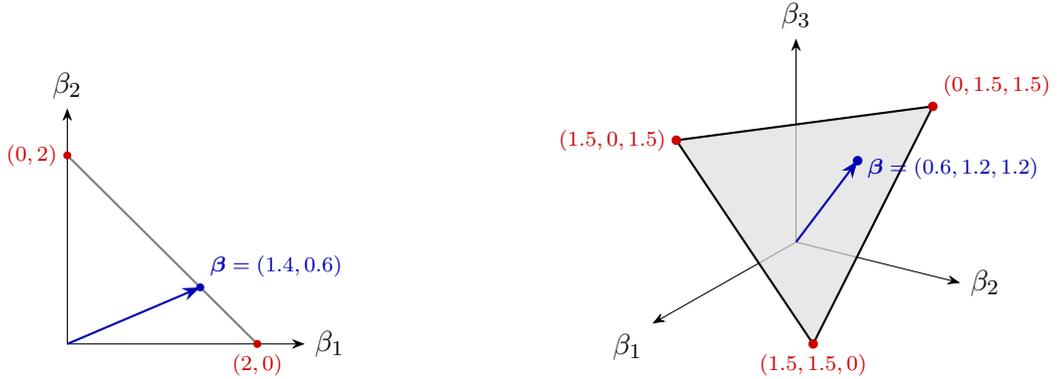
\begin{figure}[!htb]
    \vspace{0em}
    \centering
    \begin{minipage}[t]{0.48\linewidth}
        \centering
        \begin{tikzpicture}[scale=1.25, >=Stealth]
            \draw[->] (0,0) -- (2.5,0) node[right] {$\beta_1$};
            \draw[->] (0,0) -- (0,2.5) node[above] {$\beta_2$};

            \draw[thick, gray] (2,0) -- (0,2);
            \fill[red!80!black] (2,0) circle (1.2pt) node[below, font=\scriptsize, text=red!80!black] {$(2,0)$};
            \fill[red!80!black] (0,2) circle (1.2pt) node[left, font=\scriptsize, text=red!80!black] {$(0,2)$};

            \draw[->, thick, blue!70!black] (0,0) -- (1.4,0.6);
            \fill[blue!70!black] (1.4,0.6) circle (1.2pt) node[above right, font=\scriptsize] {$\bbeta=(1.4,0.6)$};
        \end{tikzpicture}
    \end{minipage}\hfill%
    \begin{minipage}[t]{0.48\linewidth}
        \centering
        \begin{tikzpicture}[
            scale=1.5,
            >=Stealth,
            x={(-0.7cm,-0.4cm)},
            y={(0.8cm,-0.2cm)},
            z={(0cm,1cm)}
        ]
        \draw[->] (0,0,0) -- (1.8,0,0) node[below left] {$\beta_1$};
        \draw[->] (0,0,0) -- (0,1.8,0) node[right] {$\beta_2$};
        \draw[->] (0,0,0) -- (0,0,1.8) node[above] {$\beta_3$};

            \coordinate (v1) at (1.5,1.5,0);
            \coordinate (v2) at (1.5,0,1.5);
            \coordinate (v3) at (0,1.5,1.5);

        \fill[gray!25, opacity=0.70] (v1)--(v2)--(v3)--cycle;
        \draw[thick] (v1)--(v2)--(v3)--cycle;

            \fill[red!80!black] (v1) circle (1.2pt) node[below, font=\scriptsize, text=red!80!black] {$(1.5,1.5,0)$};
            \fill[red!80!black] (v2) circle (1.2pt) node[left, font=\scriptsize, text=red!80!black] {$(1.5,0,1.5)$};
            \fill[red!80!black] (v3) circle (1.2pt) node[above right, font=\scriptsize, text=red!80!black] {$(0,1.5,1.5)$};

            \coordinate (b) at (0.6,1.2,1.2);
            \draw[->, thick, blue!70!black] (0,0,0) -- (b);
            \fill[blue!70!black] (b) circle (1.2pt) node[right, yshift=-3pt, font=\scriptsize] {$\bbeta=(0.6,1.2,1.2)$};
        \end{tikzpicture}
    \end{minipage}
    \caption{
        A vector $\bomega$ majorizes $\bbeta$ if and only if $\bbeta \in \conv(\mathrm{perm}(\bomega))$.
        Left ($p=2$, $k=1$): the hull is the segment between $(2,0)$ and $(0,2)$, containing $\bbeta=(1.4,0.6)$. Right ($p=3$, $k=2$): the hull is the triangle with vertices $(1.5,1.5,0)$, $(1.5,0,1.5)$, and $(0,1.5,1.5)$, containing $\bbeta=(0.6,1.2,1.2)$.
        Given $\bbeta$, Algorithm~\ref{alg:compute_g_value_root_node_algorithm} constructs such a $k$-sparse vector $\bomega$.
    }
    \label{fig:sparse_majorization_polytopes}
    \vspace{0em}
\end{figure}

The proof of Theorem~\ref{theorem:compute_g_value_algorithm_correctness} is provided in Appendix~\ref{ec_proof:compute_g_value_algorithm_correctness}.
The efficiency of Algorithm~\ref{alg:compute_g_value_root_node_algorithm} stems from exploiting the equivalence between the majorization and convex hull representations.
See Figure~\ref{fig:sparse_majorization_polytopes} for a visual illustration of this connection.
Rather than solving a generic convex optimization problem, Algorithm~\ref{alg:compute_g_value_root_node_algorithm} efficiently identifies a sparse vector $\bomega$ that majorizes $|\bar{\bbeta}|$, reducing the evaluation of a potentially complex SOCP to a simple index~search.

We now turn to the evaluation of the proximal operator. 
Calculating $\prox_{\rho g_{\mathcal{N}}}$ directly is challenging due to the implicit definition of $g_{\mathcal{N}}$. 
However, leveraging the explicit structure of the conjugate $g_{\mathcal{N}}^*$ established in Lemma~\ref{lemma:function_g_conjugate_in_BnB}, we can efficiently evaluate the proximal operator of the conjugate and recover the desired solution via the extended Moreau decomposition.

\begin{algorithm}[!tb]
    \DontPrintSemicolon
    \caption{Compute $\prox_{\rho g^*_{\mathcal{N}}}(\bbeta)$ at node $\mathcal{N}$ with associated index sets $\mathcal{J}_0, \mathcal{J}_1, \mathcal J_f$}
    \label{alg:prox_of_g_conjugate_root_node}
    \KwData{vector $\bbeta \in \mathbb R^p$, scalar $\rho > 0$, cardinality parameter $k \in [p]$, box parameter $M > 0$}
    \KwResult{$\prox_{\rho g^*_{\mathcal{N}}}(\bbeta)$}

    Set $\bar k \gets k - | \mathcal J_1 |, \ \bar p \gets | \mathcal J_f |$ \;
    Set $\bar \bbeta \gets \bbeta_{\mathcal J_f},  \ \brho \in \mathbb R_+^{\bar p}$ with $\rho_j \gets \rho$ if $j \in \{1, 2, ..., \bar k \}$ and $\rho_j \gets 0$ otherwise. \;
    Sort $\bar \bbeta$ with permutation $\bm \pi$ of $[\bar p]$ such that $\vert{\bar \beta_{\pi(1)}} \geq \ldots \geq \vert{\bar \beta_{\pi(\bar p)}} \geq 0$\;

    \lFor{$j \gets 1$ \KwTo $\bar p$}{
        $\hat{\alpha}_j \gets \prox_{\rho_j H_M} \left( \bar \beta_{\pi(j)} \right) $
    }
    $\mathcal{W} \leftarrow \{[1, 1], [2, 2]\ldots, [\bar p, \bar p]\}$. \;
    \While{$\exists [j_1, j_2], [j_2+1, j_3] \in \mathcal{W} ~ \st ~ \hat{\alpha}_{j_1} < \hat{\alpha}_{j_3}$}{
        $\mathcal{W} \gets \mathcal{W} \setminus \{[j_1, j_2]\} \setminus \{[j_2+1, j_3]\} $ \;
        $\overline{\rho} \gets \left(\sum_{j=j_1}^{j_3} \rho_j \right) / \left( j_3 - j_1 + 1 \right)$ \;
        $\overline{ \xi } \gets \left(\sum_{j=j_1}^{j_3} \bar \beta_{\pi(j)} \right) / \left( j_3 - j_1 + 1 \right)$\;
        $\hat{\alpha}_{[j_1:j_3]} \gets \prox_{\overline{\rho} H_M} \left( \overline{ \xi } \right)$\;
        $ \mathcal{W} \leftarrow \mathcal{W} \cup \{[j_1, j_3]\} $\;
    }
    Set $\balpha \in \mathbb R^p$ with $\balpha_{\mathcal{J}_0} = \bbeta_{\mathcal{J}_0}, \ \balpha_{\mathcal{J}_1} = [\prox_{\rho H_M}(\beta_j)]_{j \in \mathcal J_1}, \ \balpha_{\mathcal{J}_f} = \text{sgn}(\bar \bbeta) \odot \bpi^{-1}(\hat{\balpha})$ \;
    \KwRet{$\balpha$}
\end{algorithm}

\begin{theorem}
    \label{theorem:evaluate_prox_g_algorithm_correctness}
    For any $\bbeta \in \mathbb{R}^p$ and $\rho > 0$, Algorithm~\ref{alg:prox_of_g_conjugate_root_node} computes $\prox_{\rho g^*_{\mathcal{N}}} ( \bbeta ) $ exactly with a computational complexity of $\mathcal{O}(p + \bar p \log \bar p)$, where $\bar p = \lvert \mathcal{J}_f \rvert$ and $\bar k = k - \lvert \mathcal{J}_1 \rvert$. Moreover, Algorithm~\ref{alg:prox_of_g_conjugate_root_node} can be used as a subroutine to compute $\prox_{\rho g_{\mathcal{N}}}(\bbeta)$ exactly via the extended Moreau decomposition
    \begin{align}
        \label{eq:prox:gN}
        \prox_{\rho g_{\mathcal{N}}} ( \bbeta) = \bbeta - \rho \prox_{\rho^{-1} g_{\mathcal{N}}^*} \big( \rho^{-1} \bbeta \big).
    \end{align}
\end{theorem}

The proof of Theorem~\ref{theorem:evaluate_prox_g_algorithm_correctness} is provided in Appendix~\ref{ec_proof:evaluate_prox_g_algorithm_correctness}.
The key idea is to reformulate $\prox_{\rho g_{\mathcal{N}}^*}$ as a generalized isotonic regression problem, which minimizes a sum of Huber loss functions subject to isotonic constraints induced by the top-$\bar k$ sum operator.
This generalized isotonic regression problem is well studied in statistics and can be solved efficiently using the pool adjacent violators algorithm (PAVA) \citep{busing2022monotone}.
In Algorithm~\ref{alg:prox_of_g_conjugate_root_node}, PAVA is carried out in Lines 6--11.
In particular, adjacent violators are identified in the while loop by using an efficient block up-and-down procedure, whose detailed implementation is provided in Algorithm~\ref{alg:up_and_down_block_algorithm_for_merging_in_PAVA} in the appendix.



The computational results established in Theorems~\ref{theorem:compute_g_value_algorithm_correctness} and \ref{theorem:evaluate_prox_g_algorithm_correctness} directly satisfy the algorithmic requirements imposed by Assumption~\ref{assumption:sublinear_algorithm}.
Specifically, the efficient evaluation of the proximal operator enables the execution of standard proximal first-order methods with low per-iteration cost. 
Simultaneously, the exact computation of $g_{\mathcal{N}}(\bbeta)$ allows for the precise monitoring of the duality gap, which is the essential for the linear convergence of our proposed restart scheme.

\subsection{Convergence Guarantees for Standard GLMs}
\label{sec:conv:rate:GLM}

Having established the computational feasibility and efficiency of the algorithm, we now turn to the theoretical validity of the problem structure.
To invoke the linear convergence guarantee of Theorem~\ref{theorem:restart_linear_convergence}, we must verify that the perspective relaxation problem~\eqref{obj:original_sparse_problem_convex_composite_reformulation} satisfies the regularity conditions outlined in Assumptions~\ref{assumption:F:G}--\ref{assumption:loss_function_is_essentially_strictly_convex}.
The following proposition confirms that these assumptions hold under standard and easily verifiable conditions on the loss function.

\begin{proposition}
\label{proposition:regularity}
    Consider the primal objective $\Phi ( \bbeta ) = F ( \bX \bbeta ) + G ( \bbeta )$, where $F(\bz) := f (\bz, \by )$ with $f: \mathbb R^n \times \mathbb R^n \to \mathbb R$ being strictly convex, locally smooth, and firmly convex in $\bz$, and $G (\bbeta ) := 2 \lambda_2 \, g_{\mathcal{N}}\left(\bbeta\right)$ with $\lambda_2 > 0$ and $g_{\mathcal{N}}$ being defined in~\eqref{eq:function_g_definition}. Then, Assumptions~\ref{assumption:F:G}--\ref{assumption:loss_function_is_essentially_strictly_convex} hold.
\end{proposition}

The proof of Proposition~\ref{proposition:regularity} is provided in Appendix~\ref{ec_proof:regularity}.
We conclude with a few key remarks regarding the applicability of these results.
First, the regularity conditions established in Proposition~\ref{proposition:regularity} are satisfied by standard GLM loss functions, including least squares, logistic, and other common losses.
We provide detailed verifications in the Appendix~\ref{ec:glm_loss_verifications}.
Second, regarding the restart scheme in Theorem~\ref{theorem:restart_linear_convergence}, the analysis relies on a compact set $\mathcal{B}$ containing the iterates.
For perspective relaxation problems, the domain of the regularizer $g_{\mathcal{N}}$ is inherently compact  thanks to Lemma~\ref{lemma:g_is_proper_closed_and_convex}. 
Therefore, we can simply define $\mathcal{B} := \dom g_{\mathcal{N}}$, ensuring that the geometric constants are well-defined globally.
Finally, the exact evaluations of $G$ and its proximal, due to Theorem~\ref{theorem:compute_g_value_algorithm_correctness} and Theorem~\ref{theorem:evaluate_prox_g_algorithm_correctness}, respectively, ensure that the proposed linear convergence framework is not only theoretically sound but also computationally highly efficient for large-scale BnB applications.

\section{Experiments}
\label{sec:Experiments}

In the previous sections, we developed a principled restart scheme and low-complexity subroutines to push the computational limits of first-order methods for computing safe lower bounds.
In this section, we conduct extensive experiments on both synthetic and real-world datasets to evaluate their practical effectiveness.
We benchmark our approach against existing methods and state-of-the-art solvers to address the following empirical questions.
\begin{itemize}[label=$\diamond$,leftmargin=*]
    \item \emph{Function and Proximal Evaluation Efficiency:} How does the runtime of our specialized algorithms for evaluating $g_{\mathcal{N}}$ and $\prox_{\rho g_{\mathcal{N}}}$ compare to generic SOCP solvers?
    \item \emph{Lower Bound Computation Efficiency:} How fast can our proposed first-order method solve the perspective relaxation~\eqref{obj:perspective_relaxation} compared to state-of-the-art conic solvers? We also investigate two specific properties of our algorithmic framework.
    \begin{itemize}[label=$\bullet$,leftmargin=*]
        \item \emph{Convergence Rate:} Does our method exhibit the linear convergence rate predicted by our theory when equipped with the gap-based restart scheme?
        \item \emph{Hardware Acceleration:} Is our method GPU-friendly, and what magnitude of speedup can be achieved by leveraging GPU parallelism?
    \end{itemize}
    
    \item \emph{Optimality Certification:} When integrated into a BnB framework, how much does our method accelerate the overall process of certifying global optimality compared to existing MIP solvers?
\end{itemize}

We implemented all algorithms in Python. 
To ensure maximum efficiency, we implemented our specialized subroutines for evaluating $g_{\mathcal{N}}$ and $\prox_{\rho g_{\mathcal{N}}}$ using Numba~\citep{lam2015numba} for just-in-time compilation.
For convex baselines, we compare against the following state-of-the-art commercial and open-source SOCP solvers: Gurobi~\citep{gurobi}, MOSEK~\citep{mosek}, SCS~\citep{scs}, and Clarabel~\citep{goulart2024clarabelinteriorpointsolverconic}, using the \texttt{cvxpy} package~\citep{cvxpy} as the interface.
For the optimality certification experiments with BnB framework, we compare directly against MIP solvers Gurobi and MOSEK.

\subsection{How Fast Can We Evaluate \texorpdfstring{$g_{\mathcal{N}}(\cdot)$ and $\text{prox}_{\rho g_{\mathcal{N}}}(\cdot)$}{g and prox}?}
\begin{figure*}[!htb]
    \centering
    \includegraphics[width=1.0\textwidth]{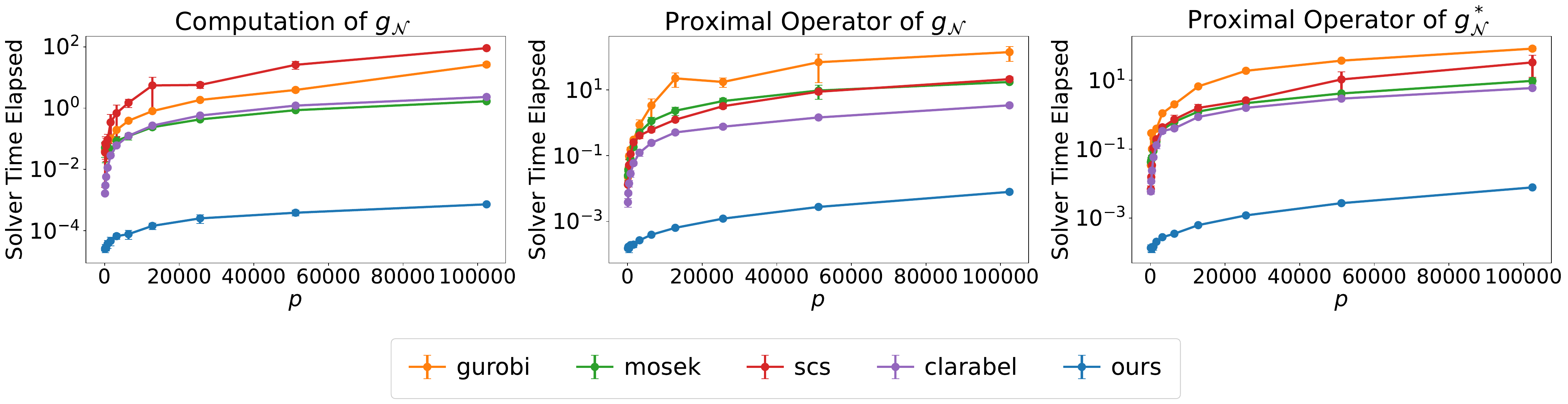}
    \caption{Running time comparison of evaluating $g_{\mathcal{N}}$ and the proximal operators of $g_{\mathcal{N}}$ and $g_{\mathcal{N}}^*$, where $\mathcal{N}$ is the root node.
    The baselines solve the corresponding SOCPs directly.}
    \label{fig:compute_g_and_prox_g_combined}
    \vspace{0em}
\end{figure*}

We benchmark our subroutines for evaluating $g_{\mathcal{N}}$ and its proximal operators against SOCP solvers that solve the corresponding SOCPs directly, where $\mathcal{N}$ is the root node.
For proximal operators, we evaluate both $\prox_{\rho g_{\mathcal{N}}} \left( \bgamma \right)$ and $\prox_{\rho^{-1} g_{\mathcal{N}}^*} \left( \rho^{-1} \bgamma \right)$.
We report running times for each method.
We repeat each setting 5 times (different random seeds) and report the mean and standard deviation (shown in error bars).
We test different methods on synthetic datasets generated as follows:
We sample the input vector $\bgamma \sim \calN(\mathbf{0}, \bI_p)$ in $\bbR^p$, where $\bI_p$ denotes the identity matrix of dimension $p$.
We vary $p \in \{2^0, \dots, 2^{10}\} \times 10^2$ and set $k=10$, $M=1.0$, and $\rho=1.0$.
When benchmarking $g_{\mathcal{N}}$, we set $\bbeta = \prox_{\rho g_{\mathcal{N}}} \left( \bgamma \right)$ by first calling Algorithm~\ref{alg:prox_of_g_conjugate_root_node}, ensuring $\bbeta \in \dom g_{\mathcal{N}}$ as required by Algorithm~\ref{alg:compute_g_value_root_node_algorithm}.

The results are shown in Figure~\ref{fig:compute_g_and_prox_g_combined}, which highlights the superiority of our method.
Compared with conventional SOCP solvers, our Algorithm~\ref{alg:compute_g_value_root_node_algorithm} achieves three orders of magnitude speedup on evaluating $g_{\mathcal{N}}$, and Algorithm~\ref{alg:prox_of_g_conjugate_root_node} together with equation~\eqref{eq:prox:gN} achieve two orders of magnitude speedup on evaluating the proximal operators. 
For instance, when $p=102400$, baselines take from several seconds to minutes to evaluate the proximal operators, whereas our approach completes the same task in $\approx 10^{-3}$ seconds for evaluating $g_{\mathcal{N}}$ and $\approx 10^{-2}$ seconds for evaluating the proximal operators, respectively.
This speed advantage, combined with the exactness of our method, will contribute to significant overall speedups when solving the perspective relaxations.
Recall that our method returns the~\textit{exact} values.
This ensures that we can calculate the primal objective function $\Phi(\bbeta)$~\textit{exactly} and apply the duality-gap-based restart scheme without any approximation errors.
Now that we can efficiently evaluate both the function value of $g_{\mathcal{N}}$ and its proximal operator, we proceed to benchmarking how fast we can compute the lower bound.

\vspace{0em}
\subsection{How Fast Can We Calculate the Lower Bound?}
\vspace{0em}

We compare against the state-of-the-art SOCP solvers for solving the perspective relaxation of the original MINLP problem, when $\mathcal{N}$ is the root node.
For our method, we use the FISTA method (with line search) equipped with our duality-gap-based restart scheme (with reduction factor $\eta=e^3$) to solve Problem~\eqref{obj:perspective_relaxation}, using our proposed algorithms to evaluate both function value of $g_{\mathcal{N}}$ and its proximal operator.
All methods are terminated upon achieving a primal-dual gap (of the BnB) tolerance of $\epsilon=10^{-6}$ or exceeding a runtime limit of 1800 seconds.
Evaluations are performed on both linear and logistic regression tasks on synthetic datasets.
For the detailed data generation process, please refer to Appendix~\ref{ec_expt:setup_for_solving_the_perspective_relaxation}.
For this experiment, we vary the feature dimension $p \in \{1000, 2000, 4000, 8000, 16000\}$.
We control the sample size by using a parameter called $n$-to-$p$ ratio, or sample to feature ratio.
We control the feature correlation by using a parameter called $\rho$.
We set the $n$-to-$p$ ratio to be $1.0$, the feature correlation $\rho$ to be $0.5$, the box constraint $M$ to be $2$, the number of nonzero coefficients (as well as the cardinality constraint) to be $10$, and the $\ell_2$ regularization coefficient $\lambda_2$ to be $1.0$.
Again, we report both the mean and standard deviation of the running time, based on 5 repeated simulations with different random seeds.

\begin{figure*}[!htb]
    \vspace{0em}
    \centering
    \includegraphics[width=0.85\textwidth]{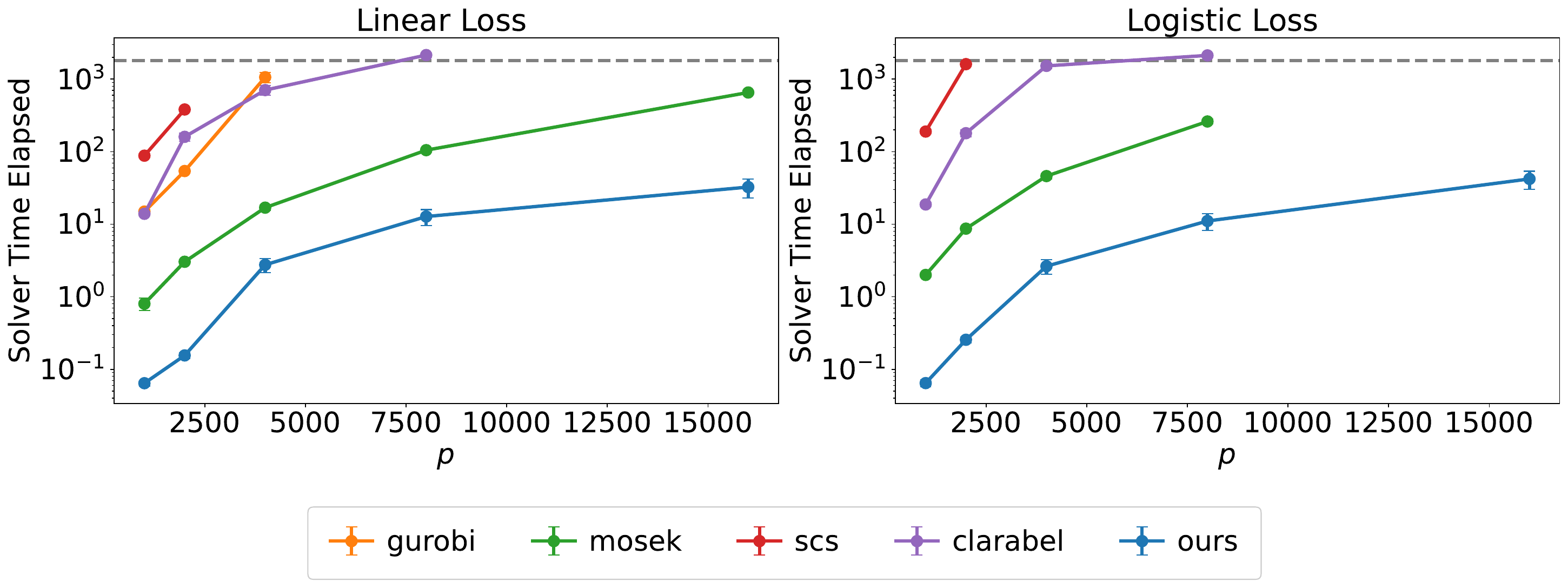}
    \caption{
        Running time comparison of solving Problem~\eqref{obj:perspective_relaxation}, the perspective relaxation at the root node.
        We set $M=2.0$, $\lambda_2=1.0$, and $n$-to-$p$ ratio to be 1.
        Gurobi cannot solve the relaxation of the cardinality constrained logistic regression problem.
    }
    \label{fig:solve_convex_relaxation_main_paper}
    \vspace{0em}
\end{figure*}

The results, shown in Figure~\ref{fig:solve_convex_relaxation_main_paper}, demonstrate that our method is faster than the fastest conic solver (MOSEK) by over one order of magnitude.
On the largest tested instances ($n=16000$ and $p=16000$), our approach attains the target tolerance ($10^{-6}$) in under 100 seconds for both the regression and classification tasks, whereas most baselines fail to converge within the 1800-second time limit.
We have also conducted additional perturbation studies, such as on the sample-to-feature ($n$-to-$p$) ratio, box constraint $M$, and $\ell_2$ regularization coefficient $\lambda_2$, which are provided in Appendix~\ref{ec_expt:perturbation_study_on_solving_the_perspective_relaxation}.
The results consistently demonstrate that our method has a significant computational advantage over existing SOCP solvers across a wide range of problem settings.
There are two factors driving this speedup.
First, both the function value and the proximal operator evaluations have low per-iteration computational costs, supported by both theory and empirical results in the previous section.
Second, our restart scheme significantly improves the convergence rate of our proximal gradient method, which we computationally quantify next.

\subsubsection{How Much Can Our Restart Scheme Accelerate Convergence?}
\begin{figure}[!h]
    \vspace{0em}
    \centering
    \includegraphics[width=1.0\textwidth]{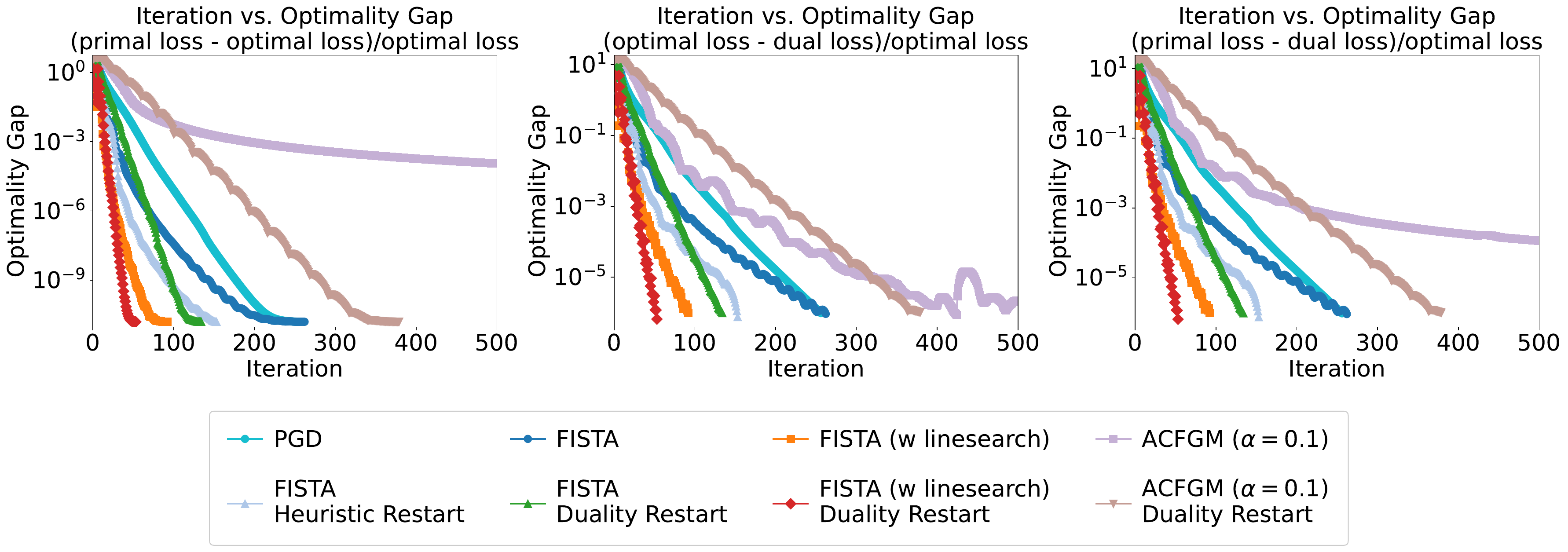}
    \caption{Convergence speed comparison between accelerated methods (with and without restart), on solving the perspective relaxation in Problem~\ref{obj:perspective_relaxation} with the linear regression loss, $n= 16000$, $p= 16000$, $k= 10$, $\rho= 0.5$, $\lambda_2 = 1.0$, and $M = 2.0$.}
    \label{fig:PDRestart_achieves_linear_convergence_rate_linear}
    \vspace{0em}
\end{figure}

\begin{figure}[!h]
    \vspace{0em}
    \centering
    \includegraphics[width=1.0\textwidth]{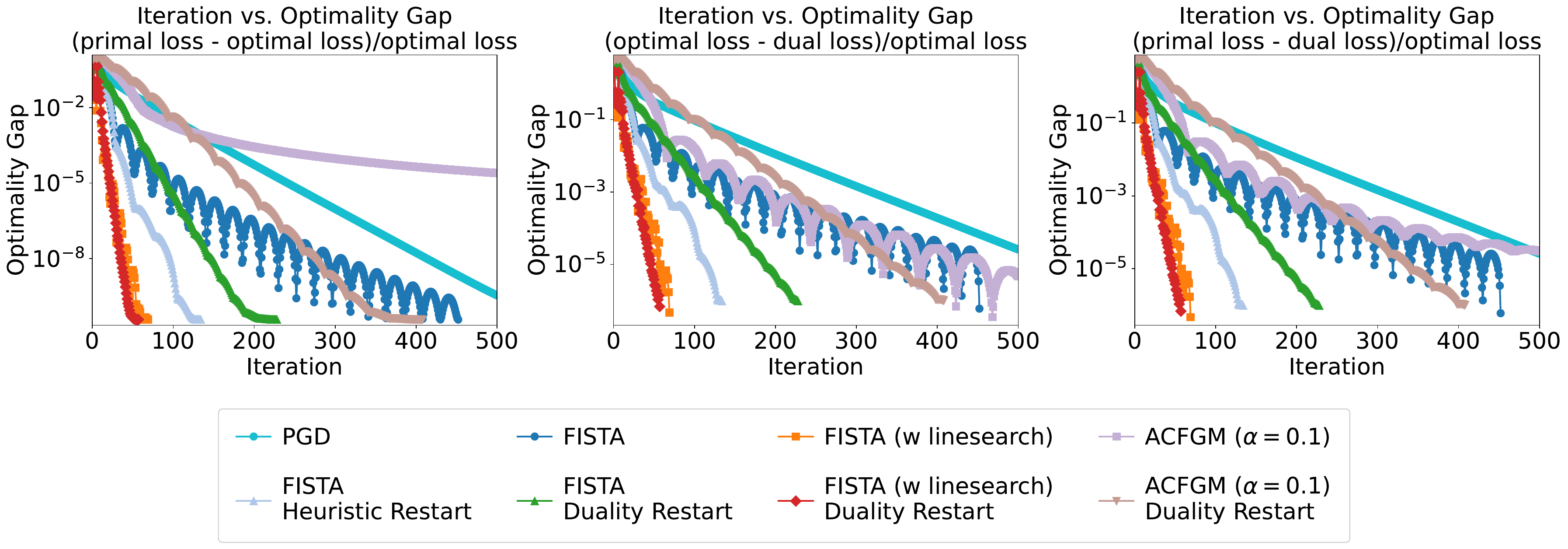}
    \caption{Convergence speed comparison between accelerated methods (with and without restart), on solving the perspective relaxation in Problem~\ref{obj:perspective_relaxation} with the logistic regression loss, $n= 16000$, $p= 16000$, $k= 10$, $\rho= 0.5$, $\lambda_2 = 1.0$, and $M = 2.0$.}
    \label{fig:PDRestart_achieves_linear_convergence_rate_logistic}
    \vspace{0em}
\end{figure}

We implemented our restart scheme on top of FISTA (both with fixed-step size and line search)~\citep{beck2017first} and ACFGM (adaptive step sizes)~\citep{li2023simple}.
For ACFGM, we implemented the algorithm as described in Corollary 2 of~\citet{li2023simple}, with $\alpha=0.1$.
We ran different methods on both the linear and logistic regression tasks, with the same parameters and data generation process as in the previous subsection.
We restarted the accelerated methods whenever the duality gap reduces by a factor of $\eta^{-1}$, where $\eta=e$.
For comparison, we also ran the vanilla versions of these methods without the restart scheme, as well as the popular proximal gradient method and FISTA with heuristic restarts based on function values~\citep{o2015adaptive}.
Figure~\ref{fig:PDRestart_achieves_linear_convergence_rate_linear} and Figure~\ref{fig:PDRestart_achieves_linear_convergence_rate_logistic} illustrate the benefits of our duality-gap-based restart scheme.
When the loss function $F$ is smooth, it is well known that the accelerated proximal gradient methods achieve a sublinear convergence rate of $O(1/T^2)$, where $T$ is the number of iterations, on the primal objective function value, which is confirmed by the results on the left plots of both figures.
When we apply our restart scheme, we observe that all accelerated methods achieve a linear convergence rate on the primal objective function value, dual objective function value (middle plot), and duality gap (right plot), matching our theoretical analysis.
Since the dual objective function value can be used as a safe lower bound, this marks, to the best of our knowledge, the first empirical demonstration of provable linear convergence for computing the safe lower bounds using a first-order method for this MINLP class.

As discussed in our theory section, our duality-gap-based restart scheme can yield linear convergence as long as $\eta>1$.
In~\ref{ec_expt:shrink_factor_in_duality_based_restart_scheme}, we provide additional numerical results showing that the practical performance of our restart scheme can be improved by increasing $\eta$ from $e$ to $e^2$ or $e^3$, in which case the duality-gap-based restart can match or outperform heuristic restart methods.
To further illustrate that our restart scheme is not limited to perspective relaxations, we also apply it to accelerate proximal gradient methods when solving the LASSO problems (both $\ell_1$-ball constrained and $\ell_1$-regularized formulations); see Appendix~\ref{ec_expt:restart_beyond_perspective_relaxations} for details.

\subsubsection{How Much Can GPUs Accelerate Lower Bound Computation?}

For fair comparison, all previous results were reported by running the experiments on CPUs.
Here, we show the potential of GPUs to further accelerate lower-bound computations.
Since the function value and proximal evaluations have low per-iteration complexity, matrix-vector multiplications dominate our proximal algorithm and can be accelerated on GPUs.
We use CuPy~\citep{cupy_learningsys2017} for GPU matrix-vector multiplications.
At each iteration, we transfer $\bgamma$ and $\bbeta$ from GPU to CPU to evaluate $\prox_{\rho g_{\mathcal{N}}} \left( \bgamma \right)$ and $g_{\mathcal{N}} \left( \bbeta \right)$ with our Numba implementation, and then transfer the results back to GPU for matrix-vector multiplications.
We rerun the experiments from the previous subsection to compare CPU-only runs with GPU-accelerated matrix-vector multiplications.
Total running times are summarized in Tables~\ref{tab:GPU_acceleration} and~\ref{tab:GPU_acceleration_logistic}.
The GPU implementation reduces runtime by an additional order of magnitude on high-dimensional instances ($p \geq 8000$) for both linear and logistic regression.

%
%
%
%
\begin{table}[!htb]
    \centering
    \vspace{0em}
    \caption{GPU acceleration of our method on the linear regression task. Top and bottom rows correspond to the mean and standard deviation of running times (seconds).}
    \vspace{0em}
    \label{tab:GPU_acceleration}
    \resizebox{0.7\linewidth}{!}{%
    \begin{tabular}{cccccc}
    \toprule
	$p$ & 1000 & 2000 & 4000 & 8000 & 16000 \\ \hline
    \multirow{2}{*}{ours CPU} & 0.070 & 0.127 & 0.447 & 5.239 & 12.491 \\
     & (0.001) & (0.002) & (0.075) & (0.632) & (1.914) \\ \hline
    \multirow{2}{*}{ours GPU} & 0.296 & 0.282 & 0.266 & 0.444 & 0.858 \\
     & (0.003) & (0.002) & (0.001) & (0.002) & (0.003)\\
     \bottomrule
    \end{tabular}%
    }
\end{table}

\begin{table}[!htb]
    \centering
    \vspace{0em}
    \caption{GPU acceleration of our method on the logistic regression task. Top and bottom rows correspond to the mean and standard deviation of running times (seconds).}
    \vspace{0em}
    \label{tab:GPU_acceleration_logistic}
    \resizebox{0.7\linewidth}{!}{%
    \begin{tabular}{cccccc}
    \toprule
	$p$ & 1000 & 2000 & 4000 & 8000 & 16000 \\ \hline
    \multirow{2}{*}{ours CPU} & 0.069 & 0.124 & 0.629 & 4.566 & 12.075 \\
     & (0.006) & (0.002) & (0.028) & (0.129) & (0.824) \\ \hline
    \multirow{2}{*}{ours GPU} & 0.275 & 0.256 & 0.249 & 0.434 & 0.997 \\
     & (0.013) & (0.015) & (0.013) & (0.017) & (0.029)\\
     \bottomrule
    \end{tabular}%
    }
    \vspace{0em}
\end{table}

\vspace{0em}
\subsection{How Fast Can Our Lower Bound Computation Certify Optimality?}

Finally, we demonstrate how our method's ability to efficiently compute safe lower bounds enables efficient optimality certification for large-scale datasets.
Note that this is just a first attempt to embed what we have into a branch and bound framework.
The goal is not to create a full-fledged MIP solver but to show the potential of using our method to compute the lower bounds to speed up the BnB process.
Integrating our lower-bound computation into a vanilla BnB framework, we prioritize node pruning via lower bound calculations while intentionally omitting advanced MIP heuristics (e.g., cutting planes, presolve routines) to evaluate the impact of our method.
We benchmark our approach against Gurobi and MOSEK, reporting runtime, optimality gaps, and the number of nodes explored.
For Gurobi on logistic regression tasks, we use the outer approximation method to approximate the logistic loss and denote this as GurobiOA.
This does not contradict Figure~\ref{fig:solve_convex_relaxation_main_paper}: for large-scale logistic relaxations, solving the root relaxation to optimality via pure cutting planes in Gurobi is computationally impractical.
For both Gurobi and MOSEK, we provide them with the same incumbent solution obtained from our BnB framework at the root node to ensure a fair comparison.
We test on both synthetic and real-world datasets, including the Santander Customer Transaction Prediction dataset~\citep{santander} for linear regression and the DOROTHEA drug discovery dataset~\citep{asuncion2007uci} for logistic regression.
For all methods, we set the optimality gap tolerance to be $10^{-6}$, the maximum runtime to be 7200 seconds, and the maximum memory to be 100 GB.
Detailed experimental configurations, including dataset descriptions, preprocessing, hyperparameter selection, and BnB implementation details, are provided in Appendix~\ref{ec_expt:setup_for_certifying_optimality_tmp}.
For the real-world datasets, we select $(\lambda_2,k,M)$ via a grid search with 5-fold cross validation using the beamsearch heuristic; see Appendix~\ref{ec_expt:cv_select_lambda2_k_M_realworld}.

In the following, we elaborate on two key aspects of using our method to calculate the safe lower bounds in the BnB framework, including a dynamic termination strategy and a warmstart technique.
Namely, we incorporate the dynamic termination strategy to avoid unnecessary proximal gradient iterations by branching or pruning nodes early.
This bypasses the need for full convergence and high-precision lower bound computations.
The detailed rules are as follows:
\begin{enumerate}
    \item If the current~\emph{primal loss} of the perspective relaxation falls below the loss of the incumbent solution ($\hat{\bbeta}$),~\textit{i.e.}, $\Phi \left( \bbeta^t \right)  \leq f \left( \bX \hat{\bbeta}, \by \right) + \lambda_2 \lVert \hat{\bbeta} \rVert_2^2$, we terminate early and proceed to the next branching step since the lower bound will never be better than the incumbent solution’s loss to prune the current node.
    \item If the current~\emph{dual loss} of the perspective relaxation exceeds the loss of the incumbent solution ($\hat{\bbeta}$),~\textit{i.e.}, $\Psi \left( \bzeta^t \right) \geq f \left( \bX \hat{\bbeta}, \by \right) + \lambda_2 \lVert \hat{\bbeta} \rVert_2^2$, we stop the lower-bound computation and prune the node immediately, since no better solution exists in its subtree.
\end{enumerate}
In addition, we also use warmstarts to speed up the lower-bound computation.
Specifically, we use the solution of the perspective relaxation problem at the parent node as the initial solution of the perspective relaxation problem at the current node.

The results are reported in Tables~\ref{tab:bnb_linear_M2}--\ref{tab:bnb_dorothea_logistic_M10}.
Our method outperforms both Gurobi and MOSEK in terms of runtime and optimality gap, often by significant margins.
In particular, on larger instances where Gurobi and MOSEK struggle to close the optimality gap within the time limit, our method consistently reaches the $0\%$ optimality gap.
On instances where all methods achieve optimality, our method often runs one to two orders of magnitude faster.

\begin{table}[!ht]
\centering
\caption{BnB results on synthetic (linear) with $k=10$, $\rho=0.5$, $\lambda_2=1.0$, $M=2.0$, $n/p=1.0$. Parentheses show lower-bound time. OOM = out of memory, TL = time $> 7200$s.}
\label{tab:bnb_linear_M2}
\resizebox{1.0\linewidth}{!}{%
\begin{tabular}{rccc ccc ccc}
\hline
 & \multicolumn{3}{c}{ours} & \multicolumn{3}{c}{gurobi} & \multicolumn{3}{c}{mosek} \\
\cmidrule(l{0.8em}r{0.8em}){2-4}\cmidrule(l{0.8em}r{0.8em}){5-7}\cmidrule(l{0.8em}r{0.8em}){8-10}
$p$ & time(s) & gap(\%) & \# nodes & time(s) & gap(\%) & \# nodes & time(s) & gap(\%) & \# nodes \\ \hline
1000 & 9 (1) & 0.00 & 105 & 146 & 0.00 & 237 & 112 & 0.00 & 113 \\
2000 & 6 (1) & 0.00 & 59 & 562 & 0.00 & 100 & 470 & 0.00 & 57 \\
4000 & 12 (6) & 0.00 & 31 & 3146 & 0.00 & 34 & 2909 & 0.00 & 31 \\
8000 & 25 (18) & 0.00 & 21 & TL & 100.00 & 1 & TL & 21.00 & 4 \\
16000 & 72 (53) & 0.00 & 21 & TL & 100.00 & 1 & OOM & OOM & OOM \\
\hline
\end{tabular}%
}
\end{table}

\begin{table}[!ht]
\centering
\caption{BnB results on synthetic (logistic) with $k=10$, $\rho=0.5$, $\lambda_2=1.0$, $M=2.0$, $n/p=1.0$. Parentheses show lower-bound time. OOM = out of memory, TL = time $> 7200$s.}
\label{tab:bnb_logistic_M2}
\resizebox{1.0\linewidth}{!}{%
\begin{tabular}{rccc ccc ccc}
\hline
 & \multicolumn{3}{c}{ours} & \multicolumn{3}{c}{gurobiOA} & \multicolumn{3}{c}{mosek} \\
\cmidrule(l{0.8em}r{0.8em}){2-4}\cmidrule(l{0.8em}r{0.8em}){5-7}\cmidrule(l{0.8em}r{0.8em}){8-10}
$p$ & time(s) & gap(\%) & \# nodes & time(s) & gap(\%) & \# nodes & time(s) & gap(\%) & \# nodes \\ \hline
1000 & 140 (7) & 0.00 & 1121 & TL & 54.17 & 597933 & 386 & 0.00 & 1083 \\
2000 & 110 (13) & 0.00 & 595 & TL & 66.26 & 582581 & 1188 & 0.00 & 583 \\
4000 & 126 (60) & 0.00 & 219 & TL & 54.04 & 255320 & 2917 & 0.00 & 235 \\
8000 & 180 (106) & 0.00 & 131 & TL & 46.95 & 81405 & OOM & OOM & OOM \\
16000 & 281 (187) & 0.00 & 67 & TL & 35.04 & 21479 & OOM & OOM & OOM \\
\hline
\end{tabular}%
}
\end{table}

\begin{table}[!ht]
\centering
\caption{BnB results on Santander (linear), $n=4459$, $p=4735$, $\lambda_2=1.0$, $M=10$. Parentheses show lower-bound time. OOM = out of memory, TL = time $> 7200$s.}
\label{tab:bnb_santander_linear_M10}
\resizebox{1.0\linewidth}{!}{%
\begin{tabular}{rccc ccc ccc}
\hline
 & \multicolumn{3}{c}{ours} & \multicolumn{3}{c}{gurobi} & \multicolumn{3}{c}{mosek} \\
\cmidrule(l{0.8em}r{0.8em}){2-4}\cmidrule(l{0.8em}r{0.8em}){5-7}\cmidrule(l{0.8em}r{0.8em}){8-10}
$k$ & time(s) & gap(\%) & \# nodes & time(s) & gap(\%) & \# nodes & time(s) & gap(\%) & \# nodes \\ \hline
6 & 345 (259) & 0.00 & 1359 & TL & 100.00 & 1 & TL & 0.33 & 437 \\
7 & 536 (395) & 0.00 & 2125 & TL & 100.00 & 1 & TL & 0.43 & 497 \\
8 & 1206 (908) & 0.00 & 4393 & TL & 100.00 & 1 & TL & 0.36 & 393 \\
9 & 3060 (2281) & 0.00 & 10085 & TL & 100.00 & 1 & TL & 0.49 & 303 \\
10 & TL (5213) & 0.01 & 18775 & TL & 100.00 & 1 & TL & 0.52 & 369 \\
\hline
\end{tabular}%
}
\end{table}

\begin{table}[!ht]
\centering
\caption{BnB results on Dorothea (logistic), $n=2300$, $p=89989$, $\lambda_2=1.0$, $M=10$. Parentheses show lower-bound time. OOM = out of memory, TL = time $> 7200$s.}
\label{tab:bnb_dorothea_logistic_M10}
\resizebox{1.0\linewidth}{!}{%
\begin{tabular}{rccc ccc ccc}
\hline
 & \multicolumn{3}{c}{ours} & \multicolumn{3}{c}{gurobiOA} & \multicolumn{3}{c}{mosek} \\
\cmidrule(l{0.8em}r{0.8em}){2-4}\cmidrule(l{0.8em}r{0.8em}){5-7}\cmidrule(l{0.8em}r{0.8em}){8-10}
$k$ & time(s) & gap(\%) & \# nodes & time(s) & gap(\%) & \# nodes & time(s) & gap(\%) & \# nodes \\ \hline
5 & 19 (15) & 0.00 & 1 & 683 & 0.00 & 938 & 1075 & 0.00 & 0 \\
10 & 66 (48) & 0.00 & 23 & 2304 & 0.00 & 3118 & OOM & OOM & OOM \\
15 & 116 (76) & 0.00 & 33 & 3005 & 0.00 & 3209 & OOM & OOM & OOM \\
20 & 223 (142) & 0.00 & 61 & TL & 0.07 & 3680 & OOM & OOM & OOM \\
25 & 629 (487) & 0.00 & 167 & TL & 0.12 & 3635 & OOM & OOM & OOM \\
\hline
\end{tabular}%
}
\end{table}

We attribute this superior performance to our method's ability to compute safe lower bounds efficiently, which is crucial for effective pruning in the BnB framework.
To support this claim, we report two additional metrics related to the BnB process: the number of nodes explored and the time spent computing lower bounds.
On instances where both our method and MOSEK achieve the $0\%$ optimality gap, our method explores comparable number of nodes but spends significantly less time computing the lower bounds, which contributes to the overall speedup.
We provide additional results on the synthetic datasets with larger box constraints ($M=3$) in~\ref{ec_expt:additional_results_on_certifying_optimality}.
The results are consistent with those when $M=2$ --- our method outperforms Gurobi and MOSEK in both runtime and optimality gap, often by significant margins.

Lastly, we demonstrate GPU acceleration for the BnB framework.
We keep the BnB implementation unchanged, except for replacing the lower-bound routine with the GPU-accelerated version from the previous subsection.
The total runtimes are summarized in Tables~\ref{tab:bnb_gpu_linear_M2}--\ref{tab:bnb_gpu_dorothea_logistic_M10}. 
For a fair comparison, CPU and GPU runs use the same machine (NVIDIA A100).
As in the root-relaxation experiments, GPU acceleration reduces lower-bound runtime by up to one order of magnitude on high-dimensional instances (Table~\ref{tab:bnb_gpu_santander_linear_M10}) for both linear and logistic tasks, which further speeds up BnB certification.
This does not always translate to an order-of-magnitude gain in total BnB time, since some work remains CPU-bound (e.g., feasible-solution search at each node); a fully GPU-native BnB implementation is beyond the scope of this work.
For smaller instances, the GPU benefit is smaller due to CPU--GPU transfer overhead.
Additional results on synthetic datasets with larger box constraints ($M=3$) are provided in~\ref{ec_expt:additional_results_on_certifying_optimality}, and are consistent with those for $M=2$.

\begin{table}[!ht]
\centering
\caption{BnB CPU vs GPU results for linear regression with $k=10$, $\rho=0.5$, $\lambda_2=1.0$, $M=2.0$, $n/p=1.0$. Parentheses show lower-bound time.}
\label{tab:bnb_gpu_linear_M2}
\resizebox{0.7\linewidth}{!}{%
\begin{tabular}{rccc ccc}
\hline
 & \multicolumn{3}{c}{ours CPU} & \multicolumn{3}{c}{ours GPU} \\
\cmidrule(l{0.8em}r{0.8em}){2-4}\cmidrule(l{0.8em}r{0.8em}){5-7}
$p$ & time(s) & gap(\%) & \# nodes & time(s) & gap(\%) & \# nodes \\ \hline
1000 & 8 (1) & 0.00 & 105 & 10 (2) & 0.00 & 105 \\
2000 & 6 (1) & 0.00 & 59 & 6 (1) & 0.00 & 59 \\
4000 & 4 (1) & 0.00 & 31 & 4 (1) & 0.00 & 31 \\
8000 & 12 (7) & 0.00 & 21 & 6 (1) & 0.00 & 21 \\
16000 & 31 (20) & 0.00 & 21 & 15 (2) & 0.00 & 21 \\
\hline
\end{tabular}%
}
\end{table}

\begin{table}[!ht]
\centering
\caption{BnB CPU vs GPU results for logistic regression with $k=10$, $\rho=0.5$, $\lambda_2=1.0$, $M=2.0$, $n/p=1.0$. Parentheses show lower-bound time.}
\label{tab:bnb_gpu_logistic_M2}
\resizebox{0.7\linewidth}{!}{%
\begin{tabular}{rccc ccc}
\hline
 & \multicolumn{3}{c}{ours CPU} & \multicolumn{3}{c}{ours GPU} \\
\cmidrule(l{0.8em}r{0.8em}){2-4}\cmidrule(l{0.8em}r{0.8em}){5-7}
$p$ & time(s) & gap(\%) & \# nodes & time(s) & gap(\%) & \# nodes \\ \hline
1000 & 132 (6) & 0.00 & 1121 & 141 (20) & 0.00 & 1123 \\
2000 & 94 (7) & 0.00 & 595 & 99 (11) & 0.00 & 595 \\
4000 & 64 (14) & 0.00 & 219 & 55 (5) & 0.00 & 219 \\
8000 & 111 (45) & 0.00 & 131 & 68 (5) & 0.00 & 131 \\
16000 & 160 (72) & 0.00 & 67 & 88 (8) & 0.00 & 67 \\
\hline
\end{tabular}%
}
\end{table}

\begin{table}[!ht]
\centering
\caption{BnB CPU vs GPU results on Santander (linear), $n=4459$, $p=4735$, $\lambda_2=1.0$, $M=10$. Parentheses show lower-bound time.}
\label{tab:bnb_gpu_santander_linear_M10}
\resizebox{0.7\linewidth}{!}{%
\begin{tabular}{rccc ccc}
\hline
 & \multicolumn{3}{c}{ours CPU} & \multicolumn{3}{c}{ours GPU} \\
\cmidrule(l{0.8em}r{0.8em}){2-4}\cmidrule(l{0.8em}r{0.8em}){5-7}
$k$ & time(s) & gap(\%) & \# nodes & time(s) & gap(\%) & \# nodes \\ \hline
6 & 169 (107) & 0.00 & 1357 & 83 (23) & 0.00 & 1363 \\
7 & 267 (166) & 0.00 & 2149 & 145 (46) & 0.00 & 2143 \\
8 & 600 (368) & 0.00 & 4385 & 337 (111) & 0.00 & 4387 \\
9 & 1655 (1058) & 0.00 & 10087 & 838 (267) & 0.00 & 10121 \\
10 & 5808 (3536) & 0.00 & 34397 & 3135 (909) & 0.00 & 34379 \\
\hline
\end{tabular}%
}
\end{table}

\begin{table}[!ht]
\centering
\caption{BnB CPU vs GPU results on Dorothea (logistic), $n=2300$, $p=89989$, $\lambda_2=1.0$, $M=10$. Parentheses show lower-bound time.}
\label{tab:bnb_gpu_dorothea_logistic_M10}
\resizebox{0.7\linewidth}{!}{%
\begin{tabular}{rccc ccc}
\hline
 & \multicolumn{3}{c}{ours CPU} & \multicolumn{3}{c}{ours GPU} \\
\cmidrule(l{0.8em}r{0.8em}){2-4}\cmidrule(l{0.8em}r{0.8em}){5-7}
$k$ & time(s) & gap(\%) & \# nodes & time(s) & gap(\%) & \# nodes \\ \hline
5 & 12 (11) & 0.00 & 1 & 3 (1) & 0.00 & 1 \\
10 & 52 (46) & 0.00 & 23 & 9 (4) & 0.00 & 23 \\
15 & 82 (69) & 0.00 & 33 & 21 (6) & 0.00 & 33 \\
20 & 167 (135) & 0.00 & 61 & 37 (11) & 0.00 & 63 \\
25 & 466 (403) & 0.00 & 167 & 84 (30) & 0.00 & 165 \\
\hline
\end{tabular}%
}
\end{table}

\section{Conclusion}
We introduce a first-order proximal framework to solve the perspective relaxations of cardinality-constrained GLM problems.
By exploiting the geometric structures of the primal and dual objective functions, we propose a new duality-gap-based restart scheme, broadly applicable beyond the specific perspective relaxation setting, to accelerate the algorithm to provably achieve the linear convergence rate, for both the primal and dual problems.
This shows, for the first time, that getting a safe lower bound can be obtained by a linearly-convergent first-order method.
Additionally, by leveraging the problem’s unique mathematical structure, we design customized subroutines to efficiently evaluate the function value and the proximal operator of the implicitly-defined regularizer $g_{\mathcal{N}}$; this opens the door to GPU acceleration since the most computationally intensive operations become matrix-vector multiplications.
Extensive empirical results demonstrate that our method outperforms state-of-the-art solvers by 1-2 orders of magnitude, establishing it as a practical, high-performance component for integration into next-generation MIP solvers.

\section*{Acknowledgements}
This work used the Delta system at the National Center for Supercomputing Applications through allocation CIS250029 from the Advanced Cyberinfrastructure Coordination Ecosystem: Services \& Support (ACCESS) program, which is supported by National Science Foundation grants \#2138259, \#2138286, \#2138307, \#2137603, and \#2138296.

\bibliographystyle{plainnat}
\bibliography{references}

@article{tillmann2024cardinality,
  title={Cardinality minimization, constraints, and regularization: a survey},
  author={Tillmann, Andreas M and Bienstock, Daniel and Lodi, Andrea and Schwartz, Alexandra},
  journal={SIAM Review},
  volume={66},
  number={3},
  pages={403--477},
  year={2024},
  publisher={SIAM}
}

@inproceedings{dash2018boolean,
  title={Boolean decision rules via column generation},
  author={Dash, Sanjeeb and G{\"u}nl{\"u}k, Oktay and Wei, Dennis},
  booktitle={Advances in Neural Information Processing Systems},
  pages={4655--4665},
  year={2018}
}

@article{lodi2024one,
  title={One-for-many counterfactual explanations by column generation},
  author={Lodi, Andrea and Ram{\'\i}rez-Ayerbe, Jasone},
  journal={arXiv:2402.09473\!\!},
  year={2024}
}

@misc{santander,
    author = {Piedra, Mercedes and Dane, Sohier and Jimenez, Soraya},
    title = {Santander Customer Transaction Prediction},
    year = {2019},
    howpublished = {\url{https://kaggle.com/competitions/santander-customer-transaction-prediction}},
    note = {{Kaggle} competition}
}

@inproceedings{cupy_learningsys2017,
  author={Okuta, Ryosuke and Unno, Yuya and Nishino, Daisuke and Hido, Shohei and Loomis, Crissman},
  title={{CuPy}: A {NumPy}-Compatible Library for {NVIDIA} {GPU} Calculations},
  booktitle={Workshop on Machine Learning Systems at Advances in Neural Information Processing Systems},
  year={2017},
}

@inproceedings{liu2022fast,
  title={Fast Sparse Classification for Generalized Linear and Additive Models},
  author={Liu, Jiachang and Zhong, Chudi and Seltzer, Margo and Rudin, Cynthia},
  booktitle={International Conference on Artificial Intelligence and Statistics},
  pages={9304--9333},
  year={2022}
}

@inproceedings{roulet2017sharpness,
  title={Sharpness, restart and acceleration},
  author={Roulet, Vincent and d'Aspremont, Alexandre},
  booktitle={Advances in Neural Information Processing Systems},
  pages={1119--1129},
  year={2017}
}

@article{fercoq2019adaptive,
  title={Adaptive restart of accelerated gradient methods under local quadratic growth condition},
  author={Fercoq, Olivier and Qu, Zheng},
  journal={IMA Journal of Numerical Analysis},
  volume={39},
  number={4},
  pages={2069--2095},
  year={2019},
  publisher={Oxford University Press}
}

@article{o2015adaptive,
  title={Adaptive restart for accelerated gradient schemes},
  author={O'Donoghue, Brendan and Candes, Emmanuel},
  journal={Foundations of Computational Mathematics},
  volume={15},
  number={3},
  pages={715--732},
  year={2015},
  publisher={Springer},
}

@article{beck2009fast,
  title={A fast iterative shrinkage-thresholding algorithm for linear inverse problems},
  author={Beck, Amir and Teboulle, Marc},
  journal={SIAM Journal on Imaging Sciences},
  volume={2},
  number={1},
  pages={183--202},
  year={2009},
  publisher={SIAM}
}

@techreport{tseng2008accelerated,
    author = {Tseng, Paul},
    title = {On accelerated proximal gradient methods for convex-concave optimization},
    year = {2008},
    institution={Department of Mathematics, University of Washington}
}

@article{nesterov1983method,
  title={A method of solving a convex programming problem with convergence rate {$O(1/k^2)$}},
  author={Nesterov, Yurii E.},
  journal={Dokl. Akad. Nauk {SSSR}},
  volume={269},
  number={3},
  pages={543--547},
  year={1983},
}

@article{chen1997convergence,
  title={Convergence rates in forward--backward splitting},
  author={Chen, George HG and Rockafellar, R Tyrrell},
  journal={SIAM Journal on Optimization},
  volume={7},
  number={2},
  pages={421--444},
  year={1997},
  publisher={SIAM}
}

@article{combettes2005signal,
  title={Signal recovery by proximal forward-backward splitting},
  author={Combettes, Patrick L and Wajs, Val{\'e}rie R},
  journal={Multiscale Modeling \& Simulation},
  volume={4},
  number={4},
  pages={1168--1200},
  year={2005},
  publisher={SIAM}
}

@inproceedings{malitsky2024adaptive,
    title={Adaptive proximal gradient method for convex optimization},
    author={Malitsky, Yura and Mishchenko, Konstantin},
    booktitle={Advances in Neural Information Processing Systems},
    pages={100670--100697},
    year={2024}
}

@article{li2023simple,
  title={A simple uniformly optimal method without line search for convex optimization},
  author={Li, Tianjiao and Lan, Guanghui},
  journal={Mathematical Programming (Forthcoming)},
  year={2025},
  publisher={Springer}
}

@inproceedings{liu2025scalable,
  title={Scalable First-order Method for Certifying Optimal k-Sparse {GLM}s},
  author={Liu, Jiachang and Shafiee, Soroosh and Lodi, Andrea},
  booktitle={Proceedings of the 42nd International Conference on Machine Learning},
  pages={39455--39481},
  year={2025}
}

@article{drusvyatskiy2018error,
  title={Error bounds, quadratic growth, and linear convergence of proximal methods},
  author={Drusvyatskiy, Dmitriy and Lewis, Adrian S},
  journal={Mathematics of Operations Research},
  volume={43},
  number={3},
  pages={919--948},
  year={2018},
  publisher={INFORMS}
}

@misc{asuncion2007uci,
  title={{The UCI Machine Learning Repository}},
  author={Asuncion, Arthur and Newman, David},
  year={2007},
  publisher={Irvine, CA, USA}
}

@inproceedings{lam2015numba,
  author={Lam, Siu Kwan and Pitrou, Antoine and Seibert, Stanley},
  editor={Finkel, Hal},
  title={Numba: a {LLVM}-based {P}ython {JIT} compiler},
  booktitle={Workshop on the {LLVM} Compiler Infrastructure in {HPC}},
  pages={7:1--7:6},
  year={2015},
}

@inproceedings{applegate2021practical,
    title={Practical large-scale linear programming using primal-dual hybrid gradient},
    author={Applegate, David and D{\'\i}az, Mateo and Hinder, Oliver and Lu, Haihao and Lubin, Miles and O'Donoghue, Brendan and Schudy, Warren},
    booktitle={Advances in Neural Information Processing Systems},
    pages={20243--20257},
    year={2021}
}

@article{lu2023cupdlp,
  title={{cuPDLP-C}: A Strengthened Implementation of {cuPDLP} for Linear Programming by {C} language},
  author={Lu, Haihao and Yang, Jinwen and Hu, Haodong and Huangfu, Qi and Liu, Jinsong and Liu, Tianhao and Ye, Yinyu and Zhang, Chuwen and Ge, Dongdong},
  year={2024},
  journal={arxiv:2312.14832\!\!}
}

@article{lu2023practical,
  title={A Practical and Optimal First-Order Method for Large-Scale Convex Quadratic Programming},
  author={Lu, Haihao and Yang, Jinwen},
  year={2025},
  journal={arxiv:2311.07710\!\!}
}

@article{han2024accelerating,
  title={Accelerating Low-Rank Factorization-Based Semidefinite Programming Algorithms on {GPU}},
  author={Han, Qiushi and Lin, Zhenwei and Liu, Hanwen and Chen, Caihua and Deng, Qi and Ge, Dongdong and Ye, Yinyu},
  year={2024},
  journal={arxiv:2407.15049\!\!},
}

@article{de2024power,
  title={On the power of linear programming for {K}-means clustering},
  author={De Rosa, Antonio and Khajavirad, Aida and Wang, Yakun},
  year={2024},
  journal={arXiv:2402.01061\!\!},
}

@book{renegar2001mathematical,
  title={A Mathematical View of Interior-Point Methods in Convex Optimization},
  author={Renegar, James},
  year={2001},
  publisher={SIAM}
}

@book{nesterov1994interior,
  title={Interior-Point Polynomial Algorithms in Convex Programming},
  author={Nesterov, Yurii and Nemirovskii, Arkadii},
  year={1994},
  publisher={SIAM}
}

@article{dikin1967iterative,
  title={Iterative solution of problems of linear and quadratic programming},
  author={Dikin, II},
  journal={Doklady Akademii Nauk},
  volume={174},
  number={4},
  pages={747--748},
  year={1967},
  organization={Russian Academy of Sciences}
}

@article{kelley1960cutting,
  title={The cutting-plane method for solving convex programs},
  author={Kelley, Jr, James E},
  journal={Journal of the Society for Industrial and Applied Mathematics},
  volume={8},
  number={4},
  pages={703--712},
  year={1960},
  publisher={SIAM}
}

@book{schrijver1998theory,
  title={Theory of Linear and Integer Programming},
  author={Schrijver, A},
  year={1998},
  publisher={John Wiley \& Sons}
}

@book{wolsey2020integer,
  title={Integer Programming},
  author={Wolsey, Laurence A},
  year={2020},
  publisher={John Wiley \& Sons}
}

@inproceedings{atamturk2020safe,
  title={Safe screening rules for {$\ell_0$}-regression from perspective relaxations},
  author={Atamturk, Alper and G{\'o}mez, Andr{\'e}s},
  booktitle={Proceedings of the 37th International Conference on Machine Learning},
  pages={421--430},
  year={2020}
}

@article{xie2020scalable,
  title={Scalable algorithms for the sparse ridge regression},
  author={Xie, Weijun and Deng, Xinwei},
  journal={SIAM Journal on Optimization},
  volume={30},
  number={4},
  pages={3359--3386},
  year={2020},
  publisher={SIAM}
}

@article{bertsimas2020sparse1,
  title={Sparse Regression: Scalable Algorithms and Empirical Performance},
  author={Bertsimas, Dimitris and Pauphilet, Jean and Van Parys, Bart},
  journal={Statistical Science},
  volume={35},
  number={4},
  pages={555--578},
  year={2020}
}

@article{bertsimas2020sparse2,
  title={Sparse high-dimensional regression: Exact scalable algorithms and phase transitions},
  author={Bertsimas, Dimitris and Van Parys, Bart},
  journal={The Annals of Statistics},
  volume={48},
  number={1},
  pages={300--323},
  year={2020}
}

@article{hazimeh2020fast,
  title={Fast best subset selection: Coordinate descent and local combinatorial optimization algorithms},
  author={Hazimeh, Hussein and Mazumder, Rahul},
  journal={Operations Research},
  volume={68},
  number={5},
  pages={1517--1537},
  year={2020},
  publisher={INFORMS}
}

@article{hazimeh2022sparse,
  title={Sparse regression at scale: Branch-and-bound rooted in first-order optimization},
  author={Hazimeh, Hussein and Mazumder, Rahul and Saab, Ali},
  journal={Mathematical Programming},
  volume={196},
  number={1},
  pages={347--388},
  year={2022},
  publisher={Springer}
}

@inproceedings{guyard2024el0ps,
  title={A New Branch-and-Bound Pruning Framework for {$\ell_0$}-Regularized Problems},
  author={Guyard, Th{\'e}o and Herzet, C{\'e}dric and Elvira, Cl{\'e}ment and Arslan, Ayse-Nur},
  booktitle={Proceedings of the 41st International Conference on Machine Learning},
  pages={48077--48096},
  year={2024}
}

@article{bertsimas2020sparse,
  title={Sparse hierarchical regression with polynomials},
  author={Bertsimas, Dimitris and Van Parys, Bart},
  journal={Machine Learning},
  volume={109},
  number={5},
  pages={973--997},
  year={2020},
  publisher={Springer}
}

@article{shafiee2024constrained,
  title={Constrained optimization of rank-one functions with indicator variables},
  author={Shafiee, Soroosh and K{\i}l{\i}n{\c{c}}-Karzan, Fatma},
  journal={Mathematical Programming},
  volume={208},
  number={1--2},
  pages={533--579},
  year={2024},
  publisher={Springer}
}

@inproceedings{zhang2023optimal,
  title={Optimal Sparse Regression Trees},
  author={Zhang, Rui and Xin, Rui and Seltzer, Margo and Rudin, Cynthia},
  booktitle={The AAAI Conference on Artificial Intelligence},
  pages={11270--11279},
  year={2023}
}

@inproceedings{liu2024fastsurvival,
  title={{FastSurvival}: Hidden Computational Blessings in Training Cox Proportional Hazards Models},
  author={Liu, Jiachang and Zhang, Rui and Rudin, Cynthia},
  booktitle={Advances in Neural Information Processing Systems},
  pages={87712--87765},
  year={2024}
}

@article{bertsimas2017optimal,
  title={Optimal classification trees},
  author={Bertsimas, Dimitris and Dunn, Jack},
  journal={Machine Learning},
  volume={106},
  number={7},
  pages={1039--1082},
  year={2017},
  publisher={Springer}
}

@inproceedings{hu2019optimal,
  title={Optimal sparse decision trees},
  author={Hu, Xiyang and Rudin, Cynthia and Seltzer, Margo},
  booktitle={Advances in Neural Information Processing Systems},
  pages={7267--7275},
  year={2019}
}

@inproceedings{liu2024okridge,
  title={O{KR}idge: Scalable optimal k-sparse ridge regression},
  author={Liu, Jiachang and Rosen, Sam and Zhong, Chudi and Rudin, Cynthia},
  booktitle={Advances in Neural Information Processing Systems},
  pages={41076--41258},
  year={2023}
}

@article{bertsimas2023learning,
  title={Learning sparse nonlinear dynamics via mixed-integer optimization},
  author={Bertsimas, Dimitris and Gurnee, Wes},
  journal={Nonlinear Dynamics},
  volume={111},
  number={7},
  pages={6585--6604},
  year={2023},
  publisher={Springer}
}

@article{ustun2016supersparse,
  title={Supersparse linear integer models for optimized medical scoring systems},
  author={Ustun, Berk and Rudin, Cynthia},
  journal={Machine Learning},
  volume={102},
  number={3},
  pages={349--391},
  year={2016},
  publisher={Springer}
}

@article{ustun2019learning,
  title={Learning Optimized Risk Scores},
  author={Ustun, Berk and Rudin, Cynthia},
  journal={Journal of Machine Learning Research},
  volume={20},
  number={150},
  pages={1--75},
  year={2019},
}

@inproceedings{liu2022fasterrisk,
  title={FasterRisk: Fast and accurate interpretable risk scores},
  author={Liu, Jiachang and Zhong, Chudi and Li, Boxuan and Seltzer, Margo and Rudin, Cynthia},
  booktitle={Advances in Neural Information Processing Systems},
  pages={17760--17773},
  year={2022}
}

@article{bienstock1996computational,
    title={Computational study of a family of mixed-integer quadratic programming problems},
    author={Bienstock, Daniel},
    journal={Mathematical Programming},
    volume={74},
    number={2},
    pages={121--140},
    year={1996},
    publisher={Springer}
}

@article{wei2022convex,
    title={On the convex hull of convex quadratic optimization problems with indicators},
    author={Wei, Linchuan and Atamt{\"u}rk, Alper and G{\'o}mez, Andr{\'e}s and K{\"u}{\c{c}}{\"u}kyavuz, Simge},
    journal={Mathematical Programming},
    volume={204},
    number={1--2},
    pages={703--737},
    year={2024}
}

@article{atamturk2020supermodularity,
    title={Supermodularity and valid inequalities for quadratic optimization with indicators},
    author={Atamt{\"u}rk, Alper and G{\'o}mez, Andr{\'e}s},
    journal={Mathematical Programming},
    volume={201},
    number={1--2},
    pages={295--338},
    year={2023},
    publisher={Springer},
}

@article{han2021compact,
    title={Compact extended formulations for low-rank functions with indicator variables},
    author={Han, Shaoning and G{\'o}mez, Andr{\'e}s},
    journal={Mathematics of Operations Research},
    volume={50},
    number={3},
    pages={1992--2016},
    year={2025},
    publisher={INFORMS}
}

@article{wei2022ideal,
    title={Ideal formulations for constrained convex optimization problems with indicator variables},
    author={Wei, Linchuan and G{\'o}mez, Andr{\'e}s and K{\"u}{\c{c}}{\"u}kyavuz, Simge},
    journal={Mathematical Programming},
    volume={192},
    number={1},
    pages={57--88},
    year={2022},
    publisher={Springer}
}

@inproceedings{wei2020convexification,
    title={On the convexification of constrained quadratic optimization problems with indicator variables},
    author={Wei, Linchuan and G{\'o}mez, Andr{\'e}s and K{\"u}{\c{c}}{\"u}kyavuz, Simge},
    booktitle={Integer Programming and Combinatorial Optimization},
    pages={433--447},
    year={2020}
}

@article{ceria1999convex,
    title={Convex programming for disjunctive convex optimization},
    author={Ceria, Sebasti{\'a}n and Soares, Jo{\~a}o},
    journal={Mathematical Programming},
    volume={86},
    number={3},
    pages={595--614},
    year={1999}
}

@article{bacci2019new,
    title={New mixed-integer nonlinear programming formulations for the unit commitment problems with ramping constraints},
    author={Bacci, Tiziano and Frangioni, Antonio and Gentile, Claudio and Tavlaridis-Gyparakis, Kostas},
    journal={Operations Research},
    volume={72},
    number={5},
    pages={2153--2167},
    year={2024},
    publisher={INFORMS}
}

@article{gunluk2010perspective,
    title={Perspective reformulations of mixed integer nonlinear programs with indicator variables},
    author={G{\"u}nl{\"u}k, Oktay and Linderoth, Jeff},
    journal={Mathematical Programming},
    volume={124},
    number={1--2},
    pages={183--205},
    year={2010},
    publisher={Springer}
}

@article{gomez2023outlier,
  title={Outlier detection in regression: {C}onic quadratic formulations},
  author={G{\'o}mez, Andr{\'e}s and Neto, Jos{\'e}},
  journal={INFORMS Journal on Computing (Forthcoming)},
  year={2025},
}

@article{gomez2021outlier,
    title={Outlier detection in time series via mixed-integer conic quadratic optimization},
    author={G{\'o}mez, Andr{\'e}s},
    journal={SIAM Journal on Optimization},
    volume={31},
    number={3},
    pages={1897--1925},
    year={2021}
}

@article{manzour2021integer,
    title={Integer programming for learning directed acyclic graphs from continuous data},
    author={Manzour, Hasan and K{\"u}{\c{c}}{\"u}kyavuz, Simge and Wu, Hao-Hsiang and Shojaie, Ali},
    journal={INFORMS Journal on Optimization},
    volume={3},
    number={1},
    pages={46--73},
    year={2021},
    publisher={INFORMS}
}

@article{dedieu2021learning,
    title={Learning Sparse Classifiers: Continuous and Mixed Integer Optimization Perspectives},
    author={Dedieu, Antoine and Hazimeh, Hussein and Mazumder, Rahul},
    journal={Journal of Machine Learning Research},
    volume={22},
    number={135},
    pages={1--47},
    year={2021},
}

@article{kucukyavuz2023consistent,
    title={Consistent Second-Order Conic Integer Programming for Learning {B}ayesian {N}etworks},
    author={Kucukyavuz, Simge and Shojaie, Ali and Manzour, Hasan and Wei, Linchuan and Wu, Hao-Hsiang},
    journal={Journal of Machine Learning Research},
    volume={24},
    number={322},
    pages={1--38},
    year={2023},
}

@article{atamturk2021sparse,
    title={Sparse and Smooth Signal Estimation: Convexification of {$\ell_0$}-Formulations},
    author={Atamt{\"u}rk, Alper and G{\'o}mez, Andr{\'e}s and Han, Shaoning},
    journal={Journal of Machine Learning Research},
    volume={22},
    number={52},
    pages={1--43},
    year={2021},
}

@book{beck2017first,
  title={First-Order Methods in Optimization},
  author={Beck, Amir},
  year={2017},
  publisher={SIAM}
}

@book{rockafellar1970convex,
    title={Convex Analysis},
    author={Rockafellar, R. Tyrrell},
    year={1970},
    publisher={Princeton University Press}
}

@article{busing2022monotone,
  title={Monotone regression: A simple and fast ${O} (n)$ {PAVA} implementation},
  author={Busing, Frank MTA},
  journal={Journal of Statistical Software},
  volume={102},
  number={Code Snippet 1},
  pages={1--25},
  year={2022}
}

@misc{scs,
  author={O'Donoghue, Brendan and Chu, Eric and Parikh, Neal and Boyd, Stephen},
  title={{SCS}: Splitting Conic Solver, version 3.2.4},
  year={2023}
}

@manual{mosek,
   author={{MOSEK ApS}},
   title={The {MOSEK} optimization toolbox for {MATLAB} manual. Version 11.0.4},
   year={2025},
}

@article{goulart2024clarabelinteriorpointsolverconic,
  title={Clarabel: An interior-point solver for conic programs with quadratic objectives},
  author={Goulart, Paul J. and Chen, Yuwen},
  year={2024},
  journal={arxiv:2405.12762\!\!},
}

@article{cvxpy,
  author={Diamond, Steven and Boyd, Stephen},
  title={{CVXPY}: {A} {Python-Embedded} Modeling Language for Convex Optimization},
  journal={Journal of Machine Learning Research},
  year={2016},
  volume={17},
  number={83},
  pages={1--5},
}

@misc{gurobi,
  author = {{Gurobi Optimization, LLC}},
  title = {{Gurobi Optimizer Reference Manual}},
  year = 2025,
}

@article{zhou2017unified,
  title={A unified approach to error bounds for structured convex optimization problems},
  author={Zhou, Zirui and So, Anthony Man-Cho},
  journal={Mathematical Programming},
  volume={165},
  number={2},
  pages={689--728},
  year={2017},
  publisher={Springer}
}

@article{luo1993error,
  title={Error bounds and convergence analysis of feasible descent methods: a general approach},
  author={Luo, Zhi-Quan and Tseng, Paul},
  journal={Annals of Operations Research},
  volume={46--47},
  number={1},
  pages={157--178},
  year={1993},
  publisher={Springer}
}

@article{beck2014fast,
  title={A fast dual proximal gradient algorithm for convex minimization and applications},
  author={Beck, Amir and Teboulle, Marc},
  journal={Operations Research Letters},
  volume={42},
  number={1},
  pages={1--6},
  year={2014},
  publisher={Elsevier}
}

@inproceedings{karimi2016linear,
  title={Linear convergence of gradient and proximal-gradient methods under the {P}olyak-{L}ojasiewicz condition},
  author={Karimi, Hamed and Nutini, Julie and Schmidt, Mark},
  booktitle={Joint European conference on machine learning and knowledge discovery in databases},
  pages={795--811},
  year={2016}
}

\input{sections/07_Appendix.tex}

\end{document}